\documentclass[10pt]{article}

\usepackage{latexsym,amsfonts}
\usepackage{amssymb,amsthm,upref,amscd}
\usepackage[T1]{fontenc}
\usepackage{times}

\newtheorem{Theorem}{Theorem}[section]
\newtheorem{definition}[Theorem]{Definition}
\newtheorem{lemma}[Theorem]{Lemma}
\newtheorem{proposition}[Theorem]{Proposition}
\newtheorem{corollary}[Theorem]{Corollary}
\newtheorem{remark}[Theorem]{Remark}
\newtheorem{Open Problem}[Theorem]{Open Problem}

\makeatletter \@addtoreset{equation}{section} \makeatother

\begin{document}

\title{\bf  Results on  entire solutions  for a degenerate critical elliptic
equation with anisotropic coefficients
 }

\author{Shaowei Chen$^{a}$  \thanks{E-mail adress:
chensw@amss.ac.cn
 } and Lishan Lin$^{b}$\\ \\
\small  $^{a}$School of Mathematical Sciences, Capital Normal
University, \\
\small Beijing 100048, P. R. China\\
\small  $^{b}$Institute of Applied Mathematics, AMSS, \\
\small Chinese Academy of Sciences, Beijing 100190, P.R. China\\
 }

\date{}
\maketitle

\begin{minipage}{13cm}
{\small {\bf Abstract:} In this paper, we study   the following
 degenerate critical elliptic equations with anisotropic
coefficients
$$
-div(|x_{N}|^{2\alpha}\nabla u)=K(x)|x_{N}|^{\alpha\cdot
2^{*}(s)-s}|u|^{2^{*}(s)-2}u\ \ \mbox{in}\ \mathbb{R}^{N}
$$
where $x=(x_{1},\cdots,x_{N})\in\mathbb{R}^{N},$ $N\geq 3,$
$\alpha>1/2,$ $0\leq s\leq 2$ and $2^{*}(s)=2(N-s)/(N-2).$ Some
basic properties of the degenerate elliptic operator
$-div(|x_{N}|^{2\alpha}\nabla u)$ are investigated and some
regularity, symmetry and uniqueness results for entire solutions
of this equation are obtained. We also get some variational
identities for solutions of this equation. As a consequence, we
obtain some nonexistence results for solutions of this equation.

}
\end{minipage}

\section{Introduction and main results}\label{diyizhang}
\hspace*{\parindent}In this paper, we study   the following
 degenerate critical elliptic equations with anisotropic
coefficients
\begin{equation}\label{yyyhsf2}
-div(|x_{N}|^{2\alpha}\nabla v)=|x_{N}|^{\alpha\cdot
2^{*}(s)-s}|v|^{2^{*}(s)-2}v\ \ \mbox{in}\ \mathbb{R}^{N}
\end{equation}
\begin{equation}\label{yyyhsf}
-div(|x_{N}|^{2\alpha}\nabla v)=K(x)|x_{N}|^{\alpha\cdot
2^{*}(s)-s}|v|^{2^{*}(s)-2}v\ \ \mbox{in}\ \mathbb{R}^{N}
\end{equation}
where $x=(x_{1},\cdots,x_{N})\in\mathbb{R}^{N},$ $N\geq 3,$
$\alpha>1/2,$ $0\leq s\leq 2$, $2^{*}(s)=2(N-s)/(N-2)$ and  $K\in
C^{1}(\mathbb{R}^{N}).$

  The motivation for studying equations (\ref{yyyhsf2}) and (\ref{yyyhsf}) comes from the
  following interesting characteristics these equations possessing.
  First, these equations  relate to the    weighted Sobolev inequality
with anisotropic coefficients  (see Theorem \ref{yyyrtrtww}):
$$\int_{\mathbb{R}^{N}}|x_{N}|^{2\alpha}|\nabla u|^{2}
\geq
C\left(\int_{\mathbb{R}^{N}}|x_{N}|^{\alpha\cdot2^{*}(s)-s}|u|^{2^{*}(s)}\right)^{2/2^{*}(s)},\
\forall u\in C^{\infty}_{0}(\mathbb{R}^{N}),$$ where $N\geq 3,$
$0\leq s\leq 2$ and $\alpha>1/2.$
 Thanks to this inequality,
 solution $u$ of  equation (\ref{yyyhsf2}) which satisfies that
$\int_{\mathbb{R}^{N}}|x_{N}|^{2\alpha}|\nabla v|^{2}<\infty,$
turns out to be a critical point of the  variational integral $J$:
$$J(v)=\frac{1}{2}\int_{\mathbb{R}^{N}}|x_{N}|^{2\alpha}|\nabla v|^{2}
-\frac{1}{2^{*}(s)}\int_{\mathbb{R}^{N}}|x_{N}|^{\alpha\cdot2^{*}(s)-s}|v|^{2^{*}(s)},\
v\in X_{\alpha}(\mathbb{R}^{N}),$$ where
$X_{\alpha}(\mathbb{R}^{N})$ is the completion space of
$C^{\infty}_{0}(\mathbb{R}^{N})$ under the the norm
$||v||=(\int_{\mathbb{R}^{N}}|x_{N}|^{2\alpha}|\nabla
v|^{2})^{1/2}$ (see Definition \ref{ttttqaa}).
   Second, equation (\ref{yyyhsf2}) is partly
conformal invariant, more precisely, if $u$ is  a solution of
equation (\ref{yyyhsf2}), then $|x|^{-(N-2+2\alpha)}u(x/|x|^{2})$
and $\mu^{(N-2+2\alpha)/2}u(\mu x+z)$ are also its solutions,
where $\mu>0$ and $z\in \mathbb{R}^{N}$ satisfying that $z_{N}=0.$
Third, these two equations are closely connected to some equations
which attracted great interest in recent years. More precisely, if
$u$ is a solution of equations (\ref{yyyhsf}), then
$v(x)=x^{\alpha}_{N}u(x),$ $x\in \mathbb{R}^{N}_{+}$ is a solution
of the following equation (see (\ref{ooo889})):
\begin{equation}\label{77yrtggfb} -\triangle
u=\frac{\lambda}{x^{2}_{N}}u+\frac{K(x)}{x^{s}_{N}}|u|^{2^{*}(s)-2}u,\
u\in \mathcal{D}^{1,2}_{0}( \mathbb{R}^{N}_{+}),\end{equation}
where $\lambda=-\alpha(\alpha-1)$ and $\mathcal{D}^{1,2}_{0}(
\mathbb{R}^{N}_{+})$ is the completion space of
$C^{\infty}_{0}(\mathbb{R}^{N})$ under the norm
$||v||=\int_{\mathbb{R}^{N}}|\nabla v|^{2}.$ Equation
(\ref{77yrtggfb}) relates to some Hardy-Sobolev inequality in half
spaces (see \cite{CL}). Let
$$H:\mathbb{R}^{N}\setminus\{(0,\cdots,0,-1)\}\rightarrow
\mathbb{R}^{N}\setminus\{(0,\cdots,0,-1)\},\ x\mapsto
\left(\frac{2x'}{1+2x_{N}+|x|^{2}},
\frac{1-|x|^{2}}{1+2x_{N}+|x|^{2}}\right)$$ and let
$\rho(x)=(2/(1+2x_{N}+|x|^{2}))^{\frac{N-2}{2}},x\in B_{1}(0).$ If
$v$ is a solution of equation (\ref{77yrtggfb}), then by
\cite{CL}, the functions $w=(u\circ H)\rho$ and
$\tilde{w}(x)=x^{\frac{N-2}{2}}_{N}u(x),$ $ x\in
\mathbb{R}^{N}_{+}$ lie in Sobolev spaces $ H^{1}_{0}( B_{1}(0))$
and $ H^{1}(\mathbb{H})$  respectively, where
 $\mathbb{H}=(\mathbb{R}^{N}_{+},dx^{2}/x^{2}_{N})$ is the
$N-$dimensional hyperbolic space,  and they are solutions of the
following two equations respectively
\begin{equation}\label{uunnvbv66tg}
-\triangle u=\frac{4\lambda}{(1-|x|^{2})^{2}} u+\frac{2^{s}K\circ
H(x)}{(1-|x|^{2})^{s}}|u|^{2^{*}(s)-2}u,\ u\in H^{1}_{0}(
B_{1}(0)),\end{equation}

\begin{equation}\label{ttfgbbcv456}
-\triangle_{\mathbb{H}^{N}}
u=\left(\lambda+\frac{N(N-2)}{4}\right)u+K(x)|u|^{2^{*}(s)-2}u,\
u\in H^{1}( \mathbb{H}).\end{equation}   In a recent paper
\cite{cfms}, the authors showed that equation (\ref{ttfgbbcv456})
can be transformed into the following equations:
\begin{itemize}
\item[(i).] semilinear elliptic equation relates to
Hardy-Sobolev-Maz'ya inequalities:
\begin{equation}\label{ij77yugggv} -\triangle
u=\frac{\mu}{|y|^{2}}u+\frac{\tilde{K}(x)}{|y|^{s}}|u|^{2^{*}(s)-2}u\
\mbox{in}\ \mathbb{R}^{N'} =\mathbb{R}^{m}\times
\mathbb{R}^{k}\end{equation} where $x=(y,z)\in\mathbb{R}^{m}\times
\mathbb{R}^{k},$ $\mu, N',$ $m, k$ depend on $N, s$, $\alpha,$ and
$\tilde{K}$ depends on $K$; \item[(ii).] Grushin type equation
with critical exponent:\begin{equation}\label{ij77yugggvuh}
-\triangle_{x}u-(\tau+1)^{2}|x|^{2\tau}\triangle_{y}u=\hat{K}(\xi)|u|^{\frac{4}{Q-2}}u\
\mbox{in}\ \mathbb{R}^{N'} =\mathbb{R}^{m}\times \mathbb{R}^{k}
\end{equation}
where $\xi=(x,y)\in\mathbb{R}^{m}\times \mathbb{R}^{k},$
$\tau,N',$ $m, k$ depend on $N, s$, $\alpha,$ and $\hat{K}$
depends on $K.$ Here $Q=m+k(1+\tau)$; \item[(iii).] semilinear
equation on Heisenberg group  and the Webster scalar curvature
equation
\begin{equation}\label{ij77y777ytrgggv}
-\triangle _{H^{N'}}u= R(\xi)|u|^{\frac{4}{Q-2}}u\ \mbox{in}\
H^{N'}
\end{equation} where $H^{N'}=
\mathbb{C}^{N'}\times\mathbb{R}=\mathbb{R}^{2N'}\times\mathbb{R},$
$Q=2N'+2,$ $\xi=(x,y,t),$ $x,y\in \mathbb{R}^{N'},$
$t\in\mathbb{R}$ and
$\triangle_{H^{N'}}=\sum^{N'}_{i=1}((\frac{\partial}{\partial
x_{i}}+2y_{i}\frac{\partial}{\partial
t})^{2}+(\frac{\partial}{\partial
y_{i}}-2x_{i}\frac{\partial}{\partial t})^{2})$. Here $N'$ depends
on $N,\alpha,s$ and $R$ depends on $K.$ \end{itemize}

A great interest has been paid to  equations
$(\ref{uunnvbv66tg})-(\ref{ij77y777ytrgggv})$ in the past years.
We refer readers  to \cite{cfms,CL, BT,CaoLi, TT,MS1,MFS,Ma,SSW}
for recent results on the existence (nonexistence), regularity,
symmetry and compactness of positive solutions of equations
$(\ref{uunnvbv66tg}) - (\ref{ij77yugggv})$. For recent development
of equations (\ref{ij77yugggvuh}) and (\ref{ij77y777ytrgggv}),
people can consult \cite{BA,FL2,FL1,M,MM} and \cite{BU,Ga,  JL,MU}
respectively.
 Through equation
(\ref{ttfgbbcv456}), equation (\ref{yyyhsf}) and equations
$(\ref{ij77yugggv})-(\ref{ij77y777ytrgggv})$  are closely linked.
Therefore,  equation (\ref{yyyhsf}) will play certain role in
studying equations $(\ref{uunnvbv66tg})-(\ref{ij77y777ytrgggv})$.

This paper is organized as follows: In section \ref{gfgbttrff}, we
obtain some weighted Sobolev type inequalities (see Theorem
\ref{yyyrtrtww}) and  define some function spaces related to these
inequalities. These inequalities can be seen as some kind of
variant of the Hardy-Sobolev-Maz'ya inequalities (see \cite{Ma}).
They not only play important role in proving the regularity and
symmetry properties of solutions of equations (\ref{yyyhsf2}) and
(\ref{yyyhsf}) but also have their own interest. In section
\ref{ggfbvt555343}, we investigate the properties of the
degenerate elliptic operator $-div(|x_{N}|^{2\alpha}\nabla u)$. We
prove a strong maximum principle (see Proposition \ref{mvxcee})
for this operator and get some results on the isolated singularity
of the positive solution of equation $-div(|x_{N}|^{2\alpha}\nabla
u)=0$ (see Proposition \ref{tttt5}). In section \ref{tgrfe},  by
means of the weighted Sobolev inequalities obtained in section 1
and the Moser iteration technique, we derive some regularity
results for positive solutions of equations (\ref{yyyhsf2}) and
(\ref{yyyhsf}). More precisely, we prove a Harnack inequality (see
Theorem \ref{nnn}) and some H\"older continuity results (see
Theorem \ref{czxzxz}) for solutions of equations (\ref{yyyhsf2})
and (\ref{yyyhsf}).  In section \ref{tgrfe},  using the moving
sphere and moving plane methods, some results on the symmetry and
decay of entire solutions of equations (\ref{yyyhsf2}) and
(\ref{yyyhsf}) are obtained (see Theorem \ref{yqwe12}, Theorem
\ref{bz2wdf} and Remark \ref{ccv5}). Especially, we obtain the
result that the positive solution of equation (\ref{yyyhsf2}) is
unique up to a M\"obius transform which leaves the upper half
space $\mathbb{R}^{N}_{+}$ invariant (see Theorem
\ref{trgt88uwhs55r}). In the last section, we derive some
variational identities (see Corollary \ref{9981}, Theorem
\ref{bvcxq} and Theorem \ref{cet551}) for solutions of equation
(\ref{yyyhsf}). As a consequence,  some non-existence results for
solutions of equation (\ref{yyyhsf}) are obtained (see Remark
\ref{ygfhcvvgf8y}).

\medskip

{\it Notation:} In what follows, $B_{\rho}(x),$
$\overline{B_{\rho}(x)}$ and $\partial B_{\rho}(x)$ will
respectively denote the open ball the closed ball and the sphere
centered at $x$ and having radius $\rho.$ For $x=(x_{1},\cdots,
x_{N})\in \mathbb{R}^{N}$, denote $(x_{1},\cdots,x_{N-1})$ by
$x'.$ The half space $\{x\in\mathbb{R}^{N}\ |\ x_{N}>0\ (<0)\}$ is
denoted by $\mathbb{R}^{N}_{+}$ (resp. $\mathbb{R}^{N}_{-}$). For
a function $u$, $u^{+}$ and $u^{-}$ denote the functions
$\max\{u(x),0\}$ and $\max\{-u(x),0\}$ respectively. For a
Lebesgue measurable set $A\subset \mathbb{R}^{N},$ $mes A$ denotes
the Lebesgue measure of $A.$  The symbol $\delta_{i,j}$ denotes
the Kronecker symbol: $\delta_{i,j}=\left\{ \matrix{1,&
  i=j\\
\cr 0, \
 & i\neq j.\cr}\right.$ For a domain $\Omega\subset
 \mathbb{R}^{N}$,
 $H^{1}_{0}(\Omega)$ is the Sobolev space defined as the completion space
 $\overline{C^{\infty}_{0}(\Omega)}^{||\cdot||}$
 under the norm $||u||=(\int_{\Omega}|u|^{2}
 +\int_{\Omega}|\nabla u|^{2})^{\frac{1}{2}}$.

\section{Some weighted  Sobolev inequalities and related function spaces }\label{gfgbttrff}
\hspace*{\parindent} In this section, we give some weighted
Sobolev type inequalities which can be seen as some kind of
variant of the Hardy-Sobolev-Maz'ya inequalities (see \cite{Ma}).
Then we define some function spaces related to these inequalities
which will be used in the subsequent sections frequently.
\begin{Theorem}\label{yyyrtrtww}
Let $N\geq 3.$ For any $0\leq s\leq 2$ and $\alpha>1/2,$ there
exist constants $C=C(\alpha,s)>0$ and $C'=C'(\alpha,s)>0$ such
that for any $u\in C^{\infty}_{0}(\mathbb{R}^{N}),$
\begin{equation}\label{yyrhf776rt}
\int_{\mathbb{R}^{N}_{\pm}}|x_{N}|^{2\alpha}|\nabla u|^{2} \geq
C'\left(\int_{\mathbb{R}^{N}_{\pm}}|x_{N}|^{\alpha\cdot2^{*}(s)-s}|u|^{2^{*}(s)}\right)^{2/2^{*}(s)};
\end{equation}
\begin{equation}\label{66rtggggvc5r}\int_{\mathbb{R}^{N}}|x_{N}|^{2\alpha}|\nabla u|^{2} \geq
C\left(\int_{\mathbb{R}^{N}}|x_{N}|^{\alpha\cdot2^{*}(s)-s}|u|^{2^{*}(s)}\right)^{2/2^{*}(s)}.
\end{equation}
\end{Theorem}
\noindent{\bf Proof.} For $u\in C^{\infty}_{0}(\mathbb{R}^{N}),$
set $v(x)=x^{\alpha}_{N}u(x),$ $x\in \mathbb{R}^{N}_{+}.$ Then
$\frac{\partial v}{\partial x_{i}}\in L^{2}(\mathbb{R}^{N}_{+})$
for $1\leq i\leq N-1$ and   $\frac{\partial v}{\partial
x_{N}}=\alpha x^{\alpha-1}_{N}u+x^{\alpha}_{N}\frac{\partial
u}{\partial x_{N}}\in L^{2}(\mathbb{R}^{N}_{+}),$ since
$\alpha>1/2.$ Therefore, $v\in H^{1}_{0}(\mathbb{R}^{N}_{+})$.
Consider
\begin{eqnarray}\label{uuuyq}
\int_{\mathbb{R}^{N}_{+}}x^{2\alpha}_{N}|\nabla u|^{2}
&=&\int_{\mathbb{R}^{N}_{+}}x^{2\alpha}_{N}|\nabla(v/x^{\alpha}_{N})|^{2}\nonumber\\
&=&\int_{\mathbb{R}^{N}_{+}}x^{2\alpha}_{N}\cdot\frac{|\nabla
v|^{2}}{x^{2\alpha}_{N}}
+\alpha^{2}\int_{\mathbb{R}^{N}_{+}}x^{2\alpha}_{N}\cdot\frac{v^{2}}{x^{2(\alpha+1)}_{N}}
-2\alpha\int_{\mathbb{R}^{N}_{+}}x^{2\alpha}_{N}\cdot\frac{v}{x^{2\alpha+1}_{N}}\cdot\frac{\partial
v}{\partial x_{N}}\nonumber\\
&=&\int_{\mathbb{R}^{N}_{+}}|\nabla
v|^{2}+\alpha^{2}\int_{\mathbb{R}^{N}_{+}}\frac{v^{2}}{x^{2}_{N}}-2\alpha\int_{\mathbb{R}^{N}_{+}}
\frac{v}{x_{N}}\cdot\frac{\partial v}{\partial x_{N}}.
\end{eqnarray}
Since $\alpha>1/2,$ we get that
$\frac{v^{2}}{x_{N}}\Big|_{x_{N}=0}=(x^{2\alpha-1}_{N}u)\Big|_{x_{N}=0}=0,$
and by the fact that the supports of $u$ and $v$  are compact, we
get that
$\frac{v^{2}}{x_{N}}\Big|_{x_{N}=\infty}=(x^{2\alpha-1}_{N}u)\Big|_{x_{N}=\infty}=0.$
Thus
\begin{eqnarray}\label{yyqqq1}
2\alpha\int_{\mathbb{R}^{N}_{+}}\frac{v}{x_{N}}\cdot\frac{\partial
v}{\partial x_{N}}
&=&2\alpha\int_{\mathbb{R}^{N-1}}\left(\int^{+\infty}_{0}\frac{v}{x_{N}}\cdot\frac{\partial
v}{\partial x_{N}}dx_{N}\right)dx'\nonumber\\
&=&\alpha\int_{\mathbb{R}^{N-1}}\left(\int^{+\infty}_{0}\frac{1}{x_{N}}\cdot\frac{\partial
}{\partial x_{N}}(v^{2})dx_{N}\right)dx'\nonumber\\
&=&\alpha\int_{\mathbb{R}^{N-1}}\left(\frac{v^{2}}{x_{N}}\Big|^{x_{N}=+\infty}_{x_{N}=0}
+\int^{+\infty}_{0}\frac{v^{2}}{x^{2}_{N}}dx_{N}\right)dx'\nonumber\\
&=&\alpha\int_{\mathbb{R}^{N}_{+}}\frac{v^{2}}{x^{2}_{N}}dx.
\end{eqnarray}
By (\ref{uuuyq}) and (\ref{yyqqq1}), we obtain
\begin{equation}\label{oooo91s}
\int_{\mathbb{R}^{N}_{+}}x^{2\alpha}_{N}|\nabla u|^{2}
=\int_{\mathbb{R}^{N}_{+}}|\nabla
v|^{2}+(\alpha^{2}-\alpha)\int_{\mathbb{R}^{N}_{+}}\frac{v^{2}}{x^{2}_{N}}.
\end{equation}
By Hardy inequality (see \cite[Theorem 327]{HL}), we have
$\int_{0}^{+\infty}\left|\frac{\partial v}{\partial
x_{N}}\right|^{2}dx_{N}\geq\frac{1}{4}\int_{0}^{+\infty}\frac{v^{2}}{x_{N}^{2}}dx_{N}.$
Since $\alpha^{2}-\alpha>-1/4,$  we get
that\begin{equation}\label{775ytgggxxc}\int_{0}^{+\infty}\left|\frac{\partial
v}{\partial
x_{N}}\right|^{2}dx_{N}+(\alpha^{2}-\alpha)\int^{+\infty}_{0}\frac{v^{2}}{x_{N}^{2}}dx_{N}
\geq \min\{1,
1+4(\alpha^{2}-\alpha)\}\int^{+\infty}_{0}\left|\frac{\partial
v}{\partial x_{N}}\right|^{2}dx_{N}.\end{equation} Notice that
$v\in H^{1}_{0}(\mathbb{R}^{N}_{+}),$  by (\ref{775ytgggxxc}) and
Hardy-Sobolev inequality in half space (see \cite{CL}), we get
that
\begin{eqnarray}\label{iiiqssx}
&&\int_{\mathbb{R}^{N}_{+}}|\nabla
v|^{2}+(\alpha^{2}-\alpha)\int_{\mathbb{R}^{N}_{+}}\frac{v^{2}}{x^{2}_{N}}
\nonumber\\
&\geq&
\int_{\mathbb{R}^{N-1}}\int_{0}^{+\infty}\left|\frac{\partial
v}{\partial
x_{N}}\right|^{2}dx_{N}dx'+(\alpha^{2}-\alpha)\int_{\mathbb{R}^{N-1}}
\int^{+\infty}_{0}\frac{v^{2}}{x_{N}^{2}}dx_{N}dx'\nonumber\\
&&+ \min\{1, 1+4(\alpha^{2}-\alpha)\}\int_{\mathbb{R}^{N-1}}
\int^{+\infty}_{0}\left|\nabla_{x'}v\right|^{2}dx_{N}dx'\nonumber\\
&\geq& \min\{1, 1+4(\alpha^{2}-\alpha)\}
\int_{\mathbb{R}^{N}_{+}}\left|\nabla v\right|^{2}dx\nonumber\\
&\geq&
C'\left(\int_{\mathbb{R}^{N}_{+}}\frac{|v|^{2^{*}(s)}}{x^{s}_{N}}\right)^{2/2^{*}(s)}=
C'\left(\int_{\mathbb{R}^{N}_{+}}|x_{N}|^{\alpha\cdot2^{*}(s)-s}|u|^{2^{*}(s)}\right)^{2/2^{*}(s)}.
\end{eqnarray}
By (\ref{iiiqssx}) and (\ref{oooo91s}), we get the inequalities
(\ref{yyrhf776rt}). The inequality (\ref{66rtggggvc5r}) follows
from  the inequalities (\ref{yyrhf776rt}) by addition.
\hfill$\Box$
\begin{definition}\label{ttttqaa} Let $\Omega$ be a bounded domain in
$\mathbb{R}^{N}$ with smooth boundary.
 Define the weighted function spaces
$X_{\alpha}(\mathbb{R}^{N})$ and $X^{0}_{\alpha}(\Omega)$ by
$$X_{\alpha}(\mathbb{R}^{N})=\overline{C^{\infty}_{0}(\mathbb{R}^{N})}
^{||\cdot||_{X_{\alpha}(\mathbb{R}^{N})}}, \
X^{0}_{\alpha}(\Omega)=\overline{C^{\infty}_{0}(\Omega)}^{||\cdot||_{X^{0}_{\alpha}(\Omega)}}$$
respectively, where the norms
$||\cdot||_{X_{\alpha}(\mathbb{R}^{N})}$ and
$||\cdot||_{X^{0}_{\alpha}(\Omega)}$ are defined by
$$||u||_{X_{\alpha}(\mathbb{R}^{N})}=(\int_{\mathbb{R}^{N}}|x_{N}|^{2\alpha}|\nabla u|^{2})^{1/2}, \
||u||_{X^{0}_{\alpha}(\Omega)}=(\int_{\Omega}|x_{N}|^{2\alpha}|\nabla
u|^{2})^{1/2}$$ for $u\in C^{\infty}_{0}(\mathbb{R}^{N})$ and
$u\in C^{\infty}_{0}(\Omega)$ respectively. By definition,
$X_{\alpha}(\mathbb{R}^{N})$ and $X^{0}_{\alpha}(\Omega)$ are
Hilbert spaces with inner products
$(u,v)=\int_{\mathbb{R}^{N}}|x_{N}|^{2\alpha}\nabla u\nabla v$ and
$(u,v)=\int_{\Omega}|x_{N}|^{2\alpha}\nabla u\nabla v$
respectively. Moreover,
 denote the space of the completion of
$C^{1}(\Omega)$ under the norm
$\left(\int_{\Omega}|x_{N}|^{2\alpha}|\nabla u|^{2} +\int_{\Omega}
|x_{N}|^{2\alpha-2} u^{2}\right)^{1/2}$ by $X_{\alpha}(\Omega)$
and denote by $X_{\alpha,loc}(\Omega)$ the space $\{u\ |\
\mbox{for any}\ D\subset\subset\Omega,\ u\in X_{\alpha}(D) \}.$
\end{definition}

\section{Some properties of  degenerate elliptic operator $-div(|x_{N}|^{2\alpha}\nabla u)$}
 \label{ggfbvt555343} \hspace*{\parindent} In this section, we investigate the
degenerate elliptic operator $-div(|x_{N}|^{2\alpha}\nabla u)$.
Throughout this section, we assume that $\alpha>1/2.$

\begin{proposition}\label{rwewe}(weak maximum principle)
If $u\in C^{2}(B_{1}(0)\setminus\{x_{N}=0\})\cap
C^{0,\gamma}(\overline{B_{1}(0)})$ for some $0<\gamma<1$ and
satisfies  \begin{equation}\label{ttegdfrrdf7y}
-div(|x_{N}|^{2\alpha}\nabla u)\geq 0\end{equation}
 weakly in $B_{1}(0),$ i.e.,
 $\int_{B_{1}(0)}|x_{N}|^{2\alpha}\nabla u\nabla\varphi\geq 0$
 for any $0\leq\varphi\in C^{\infty}_{0}(B_{1}(0))$, then
$\min_{x\in\overline{B_{1}(0)}}u(x)=\min_{x\in\partial
B_{1}(0)}u(x).$
\end{proposition}
\noindent{\bf Proof.} Without loss of generality, we may assume
that $\min_{x\in\partial B_{1}(0)}u(x)=0.$ Let $\Omega_{-}=\{x\in
B_{1}(0)\ |\ u(x)<0\}.$ If we can prove that $mes(\Omega_{-})=0,$
then the result of this Proposition holds.

Let $u^{-}(x):=\max\{-u(x),0\}$. By $\min_{x\in\partial
B_{1}(0)}u(x)=0,$ we get that $u^{-}|_{\partial B_{1}(0)}\equiv0$.
It follows that $u^{-}\in X^{0}_{\alpha}(B_{1}(0)).$ Multiplying
(\ref{ttegdfrrdf7y}) by $u^{-}$ and integrating by parts, we get
that $-\int_{B_{1}(0)}|x_{N}|^{2\alpha}|\nabla u^{-}|^{2}\geq 0.$
It follows that $u^{-}\equiv0$ in $B_{1}(0)$. Thus
$mes(\Omega_{-})=0.$ \hfill$\Box$

\medskip

Denote $e_{i}=(0,\cdots,0, \stackrel{i}{1}, 0\cdots,0),$ $1\leq
i\leq N.$

\begin{proposition}\label{mvxcee} (strong maximum principle)
Suppose that $u\in C^{2}(B_{1}(0)\setminus\{x_{N}=0\})\cap
C^{0,\gamma}(\overline{B_{1}(0)})$ for some $0<\gamma<1$. If
$-div(|x_{N}|^{2\alpha}\nabla u)\geq 0$ weakly in $B_{1}(0)$ and
$u\not\equiv constant$ in $B_{1}(0)$, then
$u(x)>\displaystyle\min_{x\in\partial B_{1}(0)}u(x),$ $x\in
B_{1}(0).$
\end{proposition}
\noindent{\bf Proof.} Without loss of generality, we may assume
that $\min_{x\in\partial B_{1}(0)}u(x)=0.$ By Proposition
  \ref{rwewe}, we know that $u\geq 0$ in $B_{1}(0).$ Since $div(|x_{N}|^{2\alpha}\nabla u)$
  is uniformly elliptic in $B_{1}(0)\setminus\{x_{N}=0\}$, by the
  classical maximum principle, we deduce that $u>0$ in
  $B_{1}(0)\setminus\{x_{N}=0\}$. Therefore, to prove this
  proposition, we only need to prove that $u(x)>0$ for $x\in
  B_{1}(0)\cap\{x_{N}=0\}.$ Without loss of generality, we only
   prove $u(0)>0.$

  Let $v(x)=|x_{N}|^{\alpha}u(x),$ $ x\in B_{1}(0).$  Straightforward calculation
shows that $$|x_{N}|^{\alpha}(\triangle v+\lambda
v/x^{2}_{N})=div(|x_{N}|^{2\alpha}\nabla u)\leq 0\ \mbox{ in}\
B_{1}(0),$$ where $\lambda=-\alpha(\alpha-1)>-1/4.$ Let
$w(x)=|x_{N}|^{\alpha}(e^{-\eta|x-a|^{2}}-e^{-\eta/4})$ with
$a=e_{N}/3$. We have
$$\triangle w+\lambda w/x^{2}_{N}=(4\eta^{2}|x-a|^{2}|x_{N}|^{\alpha}-4\eta\alpha
|x_{N}|^{\alpha}-2N\eta |x_{N}|^{\alpha}+\frac{4}{3}\eta\alpha
|x_{N}|^{\alpha-2}x_{N})e^{-\eta|x-a|^{2}}.$$ It follows that when
$\eta>0$ large enough, $\triangle w+\lambda w/x^{2}_{N}\geq0$ in
$B^+_{1/2}(a)\setminus B_{1/4}(a)$, where
$B^{+}_{1/2}(a)=B_{1/2}(a)\cap \mathbb{R}^{N}_{+}.$ Since
$\partial B_{1/4}(a)\subset B^{+}_{1}(0):=B_{1}(0)\cap
\mathbb{R}^{N}_{+}$ and $v>0$ in $B^{+}_{1}(0)$, we can choose
$\epsilon>0$ small enough such that $v(x)>\epsilon w(x)$, $x\in
\partial B_{1/4}(a).$ Thus when $\eta>0$ large enough,
$$\triangle (v-\epsilon w)+\lambda(v-\epsilon w)/x^{2}_{N}\leq 0\ \mbox{in}
\ \Omega,\quad v-w\geq 0 \ \mbox{on}\ \partial \Omega$$ where
$\Omega=B^{+}_{1/2}(a)\setminus B_{1/4}(a)$. Multiplying the above
inequality by $(v-\epsilon w)^{-}$ and integrating by parts, we
get that $-\int_{\Omega}|\nabla (v-\epsilon
w)^{-}|^{2}-\lambda\int_{\Omega}\frac{((v-\epsilon
w)^{-})^{2}}{x^{2}_{N}}\geq 0$. Since $\lambda>-1/4,$
 by Hardy inequality, we deduce that $(v-\epsilon w)^{-}=0$ in $\Omega$. Hence  $v\geq \epsilon w$ in
$\overline{\Omega}$. It follows that
$u(x)\geq\epsilon(e^{-\eta|x-a|^{2}}-e^{-\eta/4})$ for
$x\in\partial\Omega\cap\{x_{N}=0\}.$ Especially, we have $u(0)\geq
\epsilon(e^{-\eta/9}-e^{-\eta/4})>0$. \hfill$\Box$

\begin{proposition}\label{jjjhgrt}
Suppose that $u\in C^{2}(B_{1}(e_{1})\setminus\{x_{N}=0\})\cap
C^{0,\gamma}(\overline{B_{1}(e_{1})})$ for some $0<\gamma<1$. If
$-div(|x_{N}|^{2\alpha}\nabla u)\geq 0$ weakly in $B_{1}(e_{1})$,
$u(0)=\displaystyle\min_{x\in\partial B_{1}(e_{1})}u(x)=0$ and
$u>0$ in $B_{1}(e_{1}),$ then $\frac{\partial u}{\partial
x_{1}}(0)>0.$
\end{proposition}
\noindent{\bf Proof.} Let $y=x-e_{1}$ and let $v(y)=u(y+e_{1}),$
$y\in B_{1}(0).$ We have $-div(|y_{N}|^{2\alpha}\nabla
v(y))=-div(|x_{N}|^{2\alpha}\nabla u(x))\geq 0,$ weakly in
$B_{1}(0).$ Let  $w(y)=e^{-\eta |y|^{2}}-e^{-\eta}.$ We have
$\triangle w(y)=(-2N\eta+4\eta^{2}|y|^{2})e^{-\eta |y|^{2}}$ and
$\frac{\partial w}{\partial y_{N}}=-2\eta y_{N}e^{-\eta |y|^{2}}$.
Thus we get that
\begin{eqnarray}\label{uxxx}
div(|y_{N}|^{2\alpha}\nabla w) &=&|y_{N}|^{2\alpha}\triangle
w+2\alpha|y_{N}|^{2\alpha-2}y_{N}\frac{\partial w}{\partial
y_{N}}\nonumber\\
&=&(-2N\eta-4\alpha\eta+4\eta^{2}|y|^{2})|y_{N}|^{2\alpha}e^{-\eta
|y|^{2}}.\nonumber
\end{eqnarray}
When $|y|\geq 1/2$ and $\eta>0$ large enough,
$-2N\eta-4\alpha\eta+4\eta^{2}|y|^{2}\geq
-2N\eta-4\alpha\eta+\eta^{2}>0.$ Thus
$-div(|y_{N}|^{2\alpha}\nabla w)\leq 0$ in $ 1/2\leq |y|\leq 1$ if
$\eta>0$ large enough. Since  $u(y+e_{1})>0$ for any $y\in
B_{1}(0),$ we can choose $\epsilon>0$ small enough, such that
$u(y+e_{1})-\epsilon w(y)>0$ for $|y|=1/2.$ Furthermore, for
$|y|=1,$ we have $u(y+e_{1})-\epsilon w(y)\geq 0$. Thus
$$-div(|y_{N}|^{2\alpha}(u(y+e_{1})-\epsilon w(y)))\geq 0\
\mbox{in}\ 1/2\leq |y|\leq 1$$ and $$u(y+e_{1})-\epsilon w(y)\geq
0 \ \mbox{on}\ \{|y|=1/2\}\cup \{|y|=1\}.$$ By Proposition
\ref{rwewe}, we get that $u(y+e_{1})-\epsilon w(y)\geq 0 \
\mbox{in}\ 1/2\leq |y|\leq 1.$ It follows that for $0<t<1/2,$
$$\frac{u(te_{1})-u(0)}{t}=\frac{u(te_{1})}{t}
=\frac{u((t-1)e_{1}+e_{1})}{t}\geq \epsilon
\frac{w((t-1)e_{1})}{t}=\epsilon
\frac{w(te_{1}-e_{1})-w(-e_{1})}{t}.$$ Letting $t\rightarrow 0+$
in the above inequality, we get that $\frac{\partial u}{\partial
x_{1}}(0)\geq\epsilon\frac{\partial w}{\partial
y_{1}}(-e_{1})>0.$\hfill$\Box$

\bigskip

 By straightforward calculation, we
get that for any $l\in\mathbb{R},$
\begin{equation}\label{uqdd}
div(|x_{N}|^{2\alpha}\nabla(|x|^{-l}))=l(l+2-N-2\alpha)|x_{N}|^{2\alpha}|x|^{-l-2},
\ x\in \mathbb{R}^{N}\setminus \{0\}.
\end{equation}
Especially, we have
\begin{equation}\label{uq231}
div(|x_{N}|^{2\alpha}\nabla(|x|^{-(N-2+2\alpha)}))=0, \ x\in
\mathbb{R}^{N}\setminus \{0\}.
\end{equation}

For $x\in \mathbb{R}^{N},$ $r>0$,  denote $B_{r}(x)\setminus
\{x\}$ by $B^{*}_{r}(x)$.

\begin{proposition}\label{ttsd5}
Suppose that $u\in X_{\alpha,loc}(B^{*}_{2}(0))\cap
C^{0,\beta}(B^{*}_{2}(0))$ for some $\alpha>1/2$ and $0<\beta<1,$
$u>0$ in $B^{*}_{2}(0)$ and $-div(|x_{N}|^{2\alpha}\nabla u)=0$
weakly in $B^{*}_{2}(0).$ If $\lim_{|x|\rightarrow
0}|x|^{N-2+2\alpha}u(x)=0,$ then the following two results hold
\begin{itemize}\label{iaazx}
\item[(i)] there exists $M> 0$ such that $u(x)\leq M,$ $\forall
x\in B^{*}_{1}(0);$
 \item[(ii)] $u\in X_{\alpha,loc}(B_{2}(0))$
 and $-div(|x_{N}|^{2\alpha}\nabla u)=0$ weakly in $B_{2}(0).$
 \end{itemize}
\end{proposition}
\noindent{\bf Proof.} (i). Let $V_{\epsilon}(x)=\epsilon
|x|^{-(N-2+2\alpha)}+M,$ $x\in B^{*}_{1}(0)$ where $M$ is a
positive constant and $M>\sup_{x\in \partial B_{1}(0)}u(x).$ By
(\ref{uq231}), we know that $-div(|x_{N}|^{2\alpha}\nabla
V_{\epsilon})=0$ in $B^{*}_{1}(0).$ Furthermore,
$V_{\epsilon}(x)>u(x),$ $\forall x\in
\partial B_{1}(0).$ By $\lim_{|x|\rightarrow
0}|x|^{N-2+2\alpha}u(x)=0,$ we deduce that there exists a sequence
$\{\tau_{n}\}$ satisfying that $\tau_{n}\rightarrow 0+$ as
$n\rightarrow\infty$ and $V_{\epsilon}(x)>u(x),$ $\forall x\in
\partial B_{\tau_{n}}(0).$ By Proposition \ref{mvxcee}, we get
that
\begin{equation}\label{iqaz}
V_{\epsilon}(x)>u(x),\ \forall x\in B_{1}(0)\setminus
B_{\tau_{n}}(0).\end{equation}
 Fixing $n$ and letting
$\epsilon\rightarrow 0,$ by (\ref{iqaz}), we get that
\begin{equation}\label{yyyq332}
u(x)\leq M, \ \forall x\in B_{1}(0)\setminus B_{\tau_{n}}(0).
\end{equation}
Letting $n\rightarrow\infty$, by (\ref{yyyq332}), we get that
$u(x)\leq M,$ $ \forall x\in B^{*}_{1}(0).$

(ii). Let $\zeta(x)\in C^{\infty}_{0}(B_{1}(0))$ be a cut-off
function which satisfies that $0\leq \zeta\leq 1,$ in $ B_{1}(0)$,
$\zeta\equiv1$ in $B_{1/4}(0)$ and $\zeta\equiv0$ in
$\mathbb{R}^{N}\setminus B_{1/2}(0).$ Let $\eta=1-\zeta$ and
$\eta_{\epsilon}(x)=\eta(x/\epsilon).$ By
$-div(|x_{N}|^{2\alpha}\nabla u)=0$ weakly in $B^{*}_{2}(0)$, we
have $\int_{B_{1}(0)}|x_{N}|^{2\alpha}\nabla
u\nabla(\zeta\eta_{\epsilon}u)=0.$ It follows that
\begin{eqnarray}\label{urtrt}
\int_{B_{1}(0)}|x_{N}|^{2\alpha}\zeta\eta_{\epsilon} |\nabla
u|^{2}=-\int_{B_{1}(0)}|x_{N}|^{2\alpha}u\zeta\nabla
u\nabla\eta_{\epsilon} -\int_{B_{1}(0)}|x_{N}|^{2\alpha}
u\eta_{\epsilon}\nabla u\nabla\zeta.
\end{eqnarray}
We have
\begin{eqnarray}\label{0oiu}
&&\int_{B_{1}(0)}|x_{N}|^{2\alpha}u\zeta\nabla
u\nabla\eta_{\epsilon}\nonumber\\
&=&
\int_{B_{\epsilon/2}(0)\setminus
B_{\epsilon/4}(0)}|x_{N}|^{2\alpha}u\zeta\nabla
u\nabla\eta_{\epsilon}\nonumber\\
&=&\frac{1}{2}\sum^{N}_{i=1} \int_{B_{\epsilon/2}(0)\setminus
B_{\epsilon/4}(0)}|x_{N}|^{2\alpha}\zeta\frac{\partial\eta_{\epsilon}}{\partial
x_{i}}\frac{\partial }{\partial x_{i}}(u^{2})\nonumber\\
&=&\frac{1}{2}\sum^{N}_{i=1} \int_{\partial B_{\epsilon/2}(0)\cup
\partial
B_{\epsilon/4}(0)}|x_{N}|^{2\alpha}\zeta\frac{\partial\eta_{\epsilon}}{\partial
x_{i}}n_{i} \cdot u^{2}-\frac{1}{2}\sum^{N}_{i=1}
\int_{B_{\epsilon/2}(0)\setminus B_{\epsilon/4}(0)}
u^{2}\frac{\partial}{\partial
x_{i}}\left(|x_{N}|^{2\alpha}\zeta\frac{\partial\eta_{\epsilon}}{\partial
x_{i}}\right)\nonumber\\
&=& \frac{1}{2} \int_{\partial B_{\epsilon/2}(0)\cup
\partial
B_{\epsilon/4}(0)}|x_{N}|^{2\alpha}\zeta\frac{\partial\eta_{\epsilon}}{\partial
\mathbf{n}} \cdot u^{2}-\frac{1}{2}
\int_{B_{\epsilon/2}(0)\setminus B_{\epsilon/4}(0)}
u^{2}div\left(|x_{N}|^{2\alpha}\zeta\nabla\eta_{\epsilon}\right),
\end{eqnarray}
where $\mathbf{n}=(n_{1},\cdots, n_{N})$ is the outer normal
vector of $\partial B_{\epsilon/2}(0)\cup
\partial
B_{\epsilon/4}(0)$. From  result $(i)$ of this proposition,  we
know that $u$ is bounded in $B^{*}_{1}(0).$ Thus we get that
$$\lim_{\epsilon\rightarrow 0}\frac{1}{2} \int_{\partial B_{\epsilon/2}(0)\cup
\partial
B_{\epsilon/4}(0)}|x_{N}|^{2\alpha}\zeta\frac{\partial\eta_{\epsilon}}{\partial
\mathbf{n}} \cdot u^{2}=0,\ \lim_{\epsilon\rightarrow
0}\frac{1}{2} \int_{B_{\epsilon/2}(0)\setminus B_{\epsilon/4}(0)}
u^{2}div\left(|x_{N}|^{2\alpha}\zeta\nabla\eta_{\epsilon}\right)=0.$$
Thus by (\ref{0oiu}), we get \begin{equation}\label{jjjjd}
\lim_{\epsilon\rightarrow
0}\int_{B_{1}(0)}|x_{N}|^{2\alpha}u\zeta\nabla
u\nabla\eta_{\epsilon}=0.\end{equation} By (\ref{urtrt}),
(\ref{jjjjd}) and the fact that $\lim_{\epsilon\rightarrow
0}\int_{B_{1}(0)}|x_{N}|^{2\alpha} u\eta_{\epsilon}\nabla
u\nabla\zeta=\int_{B_{1}(0)}|x_{N}|^{2\alpha} u\nabla u\nabla\zeta
$,  we get that $$\int_{B_{1}(0)}|x_{N}|^{2\alpha}\zeta |\nabla
u|^{2}=\lim_{\epsilon\rightarrow
0}\int_{B_{1}(0)}|x_{N}|^{2\alpha}\zeta\eta_{\epsilon} |\nabla
u|^{2}=\int_{B_{1}(0)}|x_{N}|^{2\alpha} u\nabla
u\nabla\zeta<\infty.$$ Thus $u\in X_{\alpha}(B_{1/4}(0)).$ It
follows that $u\in X_{\alpha,loc}(B_{2}(0)).$ For any $\varphi\in
C^{\infty}_{0}(B_{2}(0)),$
 we have $$0=\int_{B_{2}(0)}|x_{N}|^{2\alpha}\nabla u\nabla (\eta_{\epsilon}\varphi)
 =\int_{B_{2}(0)}|x_{N}|^{2\alpha}\eta_{\epsilon}\nabla u\nabla\varphi
 +\int_{B_{2}(0)}|x_{N}|^{2\alpha} \varphi\nabla u\nabla\eta_{\epsilon}.$$
As the  proof of (\ref{jjjjd}), we get that
$\lim_{\epsilon\rightarrow 0}\int_{B_{2}(0)} |x_{N}|^{2\alpha}
\varphi\nabla u\nabla\eta_{\epsilon}=0.$ Moreover, we have
$$\lim_{\epsilon\rightarrow
0}\int_{B_{2}(0)}|x_{N}|^{2\alpha}\eta_{\epsilon}\nabla
u\nabla\varphi=\int_{B_{1}(0)}|x_{N}|^{2\alpha}\nabla
u\nabla\varphi.$$ Thus for any $\varphi\in
C^{\infty}_{0}(B_{2}(0)),$ $\int_{B_{2}(0)}|x_{N}|^{2\alpha}\nabla
u\nabla\varphi=0.$ \hfill$\Box$

\bigskip

The following result describes the isolated singularity of
positive solution of $-div(|x_{N}|^{2\alpha}\nabla u)=0.$  People
can consult \cite{V} for the similar result of  Laplace operator
$\triangle$.

\begin{proposition}\label{tttt5}
If $u\in X_{\alpha,loc}(B^{*}_{2}(0))\cap
C^{0,\beta}(B^{*}_{2}(0))$ for some $\alpha>1/2$ and $\beta\in
(0,1)$,  $u>0$ in $B^{*}_{2}(0)$ and $-div(|x_{N}|^{2\alpha}\nabla
u)=0$ weakly in $B^{*}_{2}(0),$ then there exists $C\geq 0$ such
that $u(x)=C|x|^{-(N-2+2\alpha)}+b(x),$ where $b(x)$ is a H\"older
continuous function in $B_{1}(0).$
\end{proposition}
\noindent{\bf Proof.} Choose $M>0$ large enough such that
$v(x)=|x|^{-(N-2+2\alpha)}-M$ satisfies $v|_{\partial
B_{1}(0)}<0.$ Let $C=\sup\{\beta\ |\ u-\beta v\geq 0\ \mbox{in}\
B^{*}_{1}(0)\}$. Obviously, $C\geq 0.$ And by the fact that there
exists $x\in B^{*}_{1}(0)$ such that $v(x)>0$, we deduce that
$C<+\infty.$

Let $w(x)=u-C v.$ For continuous function $f(x)$ defined in
$B^{*}_{1}(0)$, define $\overline{f}(r)=\max_{|x|=r}f(x),$
$\underline{f}(r)=\min_{|x|=r}f(x).$ We shall prove that
$\lim_{r\rightarrow 0}\underline{w}(r)/\overline{v}(r)=0.$

If not, there exist $\eta>0$ and $r_{n}\rightarrow 0$ such that
$\underline{w}(r_{n})\geq\eta \overline{v}(r_{n}).$ Thus $(w-\eta
v)|_{\partial B_{r_{n}}(0)}\geq 0.$ Furthermore, we have $(w-\eta
v)|_{\partial B_{1}(0)}\geq 0.$ Hence by
$-div(|x_{N}|^{2\alpha}\nabla (w-\eta v))=0$ and Proposition
\ref{rwewe}, we get that $(w-\eta v)|_{B_{1}(0)\setminus
B_{r_{n}}(0)}\geq 0.$ Letting $n\rightarrow \infty$, we get that
$(w-\eta v)|_{B^{*}_{1}(0)}\geq 0.$ It follows that
$u-(C+\eta)v\geq 0$ in $B^{*}_{1}(0)$. It contradicts the
definition of $C.$ Thus $\lim_{r\rightarrow
0}\underline{w}(r)/\overline{v}(r)=0.$ And by the Harnack
inequality (see  \cite[Theorem 4.3]{FL}), we get that
$\lim_{r\rightarrow 0}w(r)/v(r)=0.$ Then by Proposition
\ref{ttsd5} and \cite[Theorem 4.4]{FL}, we get that $w$ is a
H\"older continuous function. Let $b(x)=w(x)-CM,$ we have
$u(x)=C|x|^{-(N-2+2\alpha)}+b(x).$\hfill$\Box$

\section{Regularity   of solutions }\label{tgrfe}
\hspace*{\parindent} In this section, we derive some regularity
results for solutions of equations  (\ref{yyyhsf2}) and
(\ref{yyyhsf}).
\begin{proposition}\label{yhadsd}
If $v$ is a  solution of equation (\ref{yyyhsf2}), then
$|x|^{-(N-2+2\alpha)}v(x/|x|^{2})$ is still a solution of equation
(\ref{yyyhsf2}).
\end{proposition}
\noindent{\bf Proof.} For $x\in \mathbb{R}^{N}$, $x_{N}\geq 0,$
let $u(x)=x^{\alpha}_{N}v(x).$ By straightforward calculation, we
have $div(x^{2\alpha}_{N}\nabla v)=x^{\alpha}_{N}\triangle u
-\alpha(\alpha-1)x^{\alpha-2}_{N}u.$ Thus $u$ satisfies the
equation
\begin{equation}\label{ooo889} -\triangle
u=\frac{\lambda}{x^{2}_{N}}u+\frac{|u|^{2^{*}(s)-2}u}{x^{s}_{N}},\
x\in\mathbb{R}^{N},\ x_{N}\geq 0\end{equation} with
$\lambda=-\alpha(\alpha-1).$ Moreover, if $u$ is a
 solution of equation (\ref{ooo889}), then
$u/x^{\alpha}_{N}$ is a solution of equation (\ref{yyyhsf2}). From
\cite{CL}, we know that $|x|^{-(N-2)}u(x/|x|^{2})$ is still a
solution of equation (\ref{ooo889}). Since
$|x|^{-(N-2)}u(x/|x|^{2})=x^{\alpha}_{N}\cdot
|x|^{-(N-2+2\alpha)}v(x/|x|^{2}),$ we get that
$|x|^{-(N-2+2\alpha)}v(x/|x|^{2})$, $x\in \mathbb{R}^{N}_{+}$
satisfies  equation (\ref{yyyhsf2}). By a  similar argument, we
know that $|x|^{-(N-2+2\alpha)}v(x/|x|^{2})$, $x\in
\mathbb{R}^{N}_{-}$ also satisfies  equation (\ref{yyyhsf2}). This
finishes the proof of this proposition.\hfill$\Box$

\begin{Theorem}\label{tererd} Suppose that $\alpha>1/2$ and $0\leq s< 2.$
If $u\in X_{\alpha}(B_{1}(0))$ is a  nonnegative weak sub-solution
of equation (\ref{yyyhsf2}), i.e., for every $\varphi\in
C^{\infty}_{0}(B_{1}(0)),$ $\varphi\geq 0,$
$$\int_{B_{1}(0)}|x_{N}|^{2\alpha}\nabla u\nabla\varphi\leq
\int_{B_{1}(0)}|x_{N}|^{\alpha\cdot
2^{*}(s)-s}u^{2^{*}(s)-1}\varphi,$$ then there exists $\sigma\in
(0,1)$ such that $u\in L^{\infty}(B_{\sigma}(0)).$
\end{Theorem}
\noindent{\bf Proof.} For $t>2,\ k>0,$ define $h(r)=\left\{
\matrix{ r^{t/2},&
 \ 0\leq
r\leq k\\
\cr {t\over{2}}k^{{t\over{2}}-1}r+(1-{t\over{2}})k^{t\over{2}}, \
 & r\geq k\\\cr}\right. $, $\phi(r)=\int^{r}_{0}|h'(s)|^{2}ds.$
It is easy to verify the following two inequalities
\begin{equation}\label{3c}
|r\phi(r)|\leq {{t^{2}}\over{4(t-1)}}|h(r)|^{2},
\end{equation}
 \begin{equation} \label{3d} |\phi(r)-h(r)h'(r)|\leq
C_{t}|h(r)h'(r)|,
\end{equation}
where $C_{t}={{t-2}\over{2(t-1)}}<1.$
 Let $0<\tau<\rho<1$. Choose
$\eta \in C^{\infty}_{0}(B_{\rho}(0))$ satisfying $ 0\leq\eta\leq
1,$ $ \eta\equiv 1 $ in $ B_{\tau}(0),$ $ \eta\equiv 0$ in $
\mathbb{R}^{N}\setminus B_{\rho}(0)$ and $|\nabla \eta|\leq
2/(\rho-\tau). $ Then $\eta^{2}\phi(u),\eta h(u) \in
X_{\alpha}^{0}(B_{1}(0)).$ We have
\begin{eqnarray}
\int_{B_{1}(0)}|x_{N}|^{2\alpha}\nabla u\nabla (\eta^{2}\phi(u))
&=&\int_{B_{1}(0)}|x_{N}|^{2\alpha}\eta^{2}(h'(u))^{2}|\nabla u
|^{2}+2\int_{B_{1}(0)}|x_{N}|^{2\alpha}\eta\phi(u)\nabla u\nabla
\eta\nonumber\\
&=&\int_{B_{1}(0)}|x_{N}|^{2\alpha}\eta^{2}|\nabla ( h(u))|^{2}+
2\int_{B_{1}(0)}|x_{N}|^{2\alpha}\eta\phi(u)\nabla u\nabla
\eta.\nonumber
\end{eqnarray}
Notice that $|\nabla(\eta h(u))|^{2}=\eta^{2}|\nabla
(h(u))|^{2}+h^{2}(u)|\nabla\eta|^{2}+2\eta h(u)\nabla
(h(u))\nabla\eta,$ by (\ref{3d}), we have
\begin{eqnarray}\label{3i}
\int_{B_{1}(0)}|x_{N}|^{2\alpha}\nabla u\nabla
(\eta^{2}\phi(u))&=& \int_{B_{1}(0)}|x_{N}|^{2\alpha}|\nabla (\eta
h(u))|^{2}
-\int_{B_{1}(0)}|x_{N}|^{2\alpha}h^{2}(u)|\nabla \eta|^{2}\nonumber\\
&&- 2\int_{B_{1}(0)}|x_{N}|^{2\alpha}\eta h(u)h'(u)\nabla
u\nabla\eta +
2\int_{B_{1}(0)}|x_{N}|^{2\alpha}\eta\phi(u)\nabla u\nabla\eta\nonumber\\
&\geq&\int_{B_{1}(0)}|x_{N}|^{2\alpha}|\nabla (\eta h(u))|^{2}-
\int_{B_{1}(0)}|x_{N}|^{2\alpha}h^{2}(u)|\nabla \eta|^{2}\nonumber\\
&&-
2\int_{B_{1}(0)}|x_{N}|^{2\alpha}\eta|\phi(u)-h(u)h'(u)|\cdot|\nabla
u\nabla\eta|\nonumber\\
&\geq& \int_{B_{1}(0)}|x_{N}|^{2\alpha}|\nabla(\eta h(u))|^{2}-
\int_{B_{1}(0)}|x_{N}|^{2\alpha}h^{2}(u)|\nabla \eta|^{2} \nonumber\\
&&- 2C_{t}\int_{B_{1}(0)}|x_{N}|^{2\alpha}|\eta
h(u)\nabla(h(u))\nabla \eta|.
\end{eqnarray}
Since
\begin{eqnarray}\label{zzz}
\int_{B_{1}(0)}|x_{N}|^{2\alpha}|\eta h(u)\nabla (h(u))\nabla
\eta|
 &=& \int_{B_{1}(0)}|x_{N}|^{2\alpha}|(\nabla (\eta
h(u))-h(u)\nabla\eta)\nabla\eta|\cdot|h(u)|\nonumber\\
&\leq& \int_{B_{1}(0)}|x_{N}|^{2\alpha}|h(u)\nabla(\eta
h(u))\nabla\eta|+
\int_{B_{1}(0)}|x|^{2\alpha}|h(u)|^{2}|\nabla \eta|^{2}\nonumber\\
&\leq&{1\over{2}}\int_{B_{1}(0)}|x_{N}|^{2\alpha}h^{2}(u)|\nabla\eta|^{2}
+ {1\over{2}}\int_{B_{1}(0)}|x_{N}|^{2\alpha}|\nabla(\eta
h(u))|^{2}
\nonumber\\
&&+\int_{B_{1}(0)}|x_{N}|^{2\alpha}|h(u)|^{2}|\nabla\eta|^{2},
\end{eqnarray}
 by (\ref{3i}),  (\ref{zzz}) and the weighted inequality (\ref{66rtggggvc5r}), we deduce
 that
\begin{eqnarray}\label{3f}
&&\int_{B_{1}(0)}|x_{N}|^{2\alpha}\nabla u\nabla (\eta^{2}\phi(u))\nonumber\\
&\geq& \int_{B_{1}(0)}|x_{N}|^{2\alpha}|\nabla (\eta
h(u))|^{2}-\int_{B_{1}(0)}|x_{N}|^{2\alpha}h^{2}(u)|\nabla\eta|^{2}\nonumber\\
&&-2C_{t}
\left({3\over{2}}\int_{B_{1}(0)}|x_{N}|^{2\alpha}h^{2}(u)|\nabla\eta|^{2}
+ {1\over{2}}\int_{B_{1}(0)}|x_{N}|^{2\alpha}|\nabla(\eta
h(u))|^{2} \right)\nonumber\\
&=&{t\over{2(t-1)}}\int_{B_{1}(0)}|x_{N}|^{2\alpha}|\nabla (\eta
h(u))|^{2}-(1+3C_{t})\int_{B_{1}(0)}|x_{N}|^{2\alpha}h^{2}(u)|\nabla
\eta|^{2}\nonumber\\
&\geq&{Ct\over{2(t-1)}}(\int_{B_{1}(0)}|x_{N}|^{\alpha\cdot
2^{*}(s)-s}|\eta h(u)|^{2^{*}(s)})^{\frac{2}{2^{*}(s)}}
-(1+3C_{t})\int_{B_{1}(0)}|x_{N}|^{2\alpha}h^{2}(u)|\nabla
\eta|^{2}.
\end{eqnarray}
By (\ref{3c}) and  H\"older inequality, we have
\begin{eqnarray}\label{3e}
&&\int_{B_{1}(0)}|x_{N}|^{\alpha\cdot
2^{*}(s)-s}u^{2^{*}(s)-1}\eta^{2}\phi(u)\nonumber\\ &\leq&
{{t^{2}}\over{4(t-1)}}\int_{B_{1}(0)}|x_{N}|^{\alpha\cdot
2^{*}(s)-s}|u|^{2^{*}(s)-2}|\eta
h(u)|^{2}\nonumber\\
&\leq& {{t^{2}}\over{4(t-1)}}\left(\int_{\eta\neq
0}|x_{N}|^{\alpha\cdot
2^{*}(s)-s}|u|^{2^{*}(s)}\right)^{\frac{2^{*}(s)-2}{2^{*}(s)}}
\left(\int_{B_{1}(0)}|x_{N}|^{\alpha\cdot 2^{*}(s)-s}|\eta h(
u)|^{2^{*}(s)}\right)^{\frac{2}{2^{*}(s)}}
\end{eqnarray}
Since $u$ is a nonnegative weak sub-solution of equation
(\ref{yyyhsf2}), we have
$$\int_{B_{1}(0)}|x_{N}|^{2\alpha}\nabla u\nabla
(\eta^{2}\phi(u))\leq\int_{B_{1}(0)}|x_{N}|^{\alpha\cdot
2^{*}(s)-s}u^{2^{*}(s)-1}\eta^{2}\phi(u).$$ Then by (\ref{3f}) and
(\ref{3e}) we get that
\begin{eqnarray}\label{3g}
&&\left(\int_{B_{1}(0)}|x_{N}|^{\alpha\cdot 2^{*}(s)-s}|\eta
h(u)|^{2^{*}(s)}\right)^{2\over{2^{*}(s)}}\nonumber\\
&\leq& {{t}\over{2C}} \left(\int_{\eta\neq 0}|x_{N}|^{\alpha\cdot
2^{*}(s)-s}|u|^{2^{*}(s)}\right)^{{2^{*}(s)-2}\over{2^{*}(s)}}
\left(\int_{B_{1}(0)}|x_{N}|^{\alpha\cdot 2^{*}(s)-s}|\eta h(
u)|^{2^{*}(s)}\right)^{2\over{2^{*}(s)}}\nonumber\\
&&+{{2(1+3C_{t})(t-1)}\over{Ct}}\int_{B_{1}(0)}
|x_{N}|^{2\alpha}h^{2}(u)|\nabla \eta |^{2}.
\end{eqnarray}
Choose $\rho$ small enough such that ${{t}\over{2C}}
\left(\int_{\eta\neq 0}|x_{N}|^{\alpha\cdot
2^{*}(s)-s}|u|^{2^{*}(s)}\right)^{{2^{*}(s)-2}\over{2^{*}(s)}}<1/2.$
Notice that ${{2(1+3C_{t})(t-1)}\over{t}}<8$ (since $0<C_{t}<1$
and $t>2$) and $|\nabla\eta|<2/(\rho-\tau),$ from (\ref{3g}) we
have
\begin{equation}\label{3h} \left(\int_{B_{\tau}(0)}|x_{N}|^{\alpha\cdot2^{*}(s)-s}|
h(u)|^{2^{*}(s)}\right)^{2\over{2^{*}(s)}} \leq
{{64}\over{C(\rho-\tau)^{2}}}
\int_{B_{\rho}(0)}|x_{N}|^{2\alpha}h^{2}(u).
\end{equation}
Choose $t_{0}>2$ such that $t_{0}-2$ small enough and let
$k\rightarrow\infty$ in (\ref{3h}), we get
\begin{equation}\label{rt}
\left(\int_{B_{\tau}(0)}|x_{N}|^{\alpha\cdot 2^{*}(s)-s}
|u|^{2^{*}(s)t_{0}/2}\right)^{2\over{2^{*}(s)}} \leq
{{64}\over{C(\rho-\tau)^{2}}}
\int_{B_{\rho}(0)}|x_{N}|^{2\alpha}|u|^{t_{0}}. \end{equation} Let
$s_{0}\in (0,2)$ be such that $2^{*}(s_{0})=t_{0}.$ Then
$s_{0}\rightarrow 2$ as $t_{0}\rightarrow 2.$ It follows that
$2\alpha>\alpha\cdot 2^{*}(s_{0})-s_{0}$ if $t_{0}-2>0$ small
enough. Thus by Theorem \ref{yyyrtrtww}, we get that
$(\int_{B_{1}(0)}|x_{N}|^{2\alpha}|\zeta
u|^{t_{0}})^{\frac{2}{t_{0}}}\leq(\int_{B_{1}(0)}|x_{N}|^{\alpha\cdot
2^{*}(s_{0})-s_{0}}|\zeta u|^{t_{0}})^{\frac{2}{t_{0}}}\leq
C\int_{B_{1}(0)}|x_{N}|^{2\alpha}|\nabla (\zeta u)|^{2}<\infty,$
where $\zeta\in C^{\infty}_{0}(B_{1}(0))$ is a cut-off function
with $\zeta\equiv 1$ in $B_{\rho}(0)$. Combining (\ref{rt}), we
get that
\begin{equation}\label{yy} \int_{B_{\rho}(0)}|x_{N}|^{\alpha\cdot2^{*}(s)-s}
|u|^{2^{*}(s)t_{0}/2}<\infty.
\end{equation}

For any $0<r_{2}<r_{1}\leq\rho,$ let $\eta\in
C^{\infty}_{0}(B_{r_{1}})$ be a cut-off function which satisfies
that $ 0\leq\eta\leq 1,$ $ \eta\equiv 1 $ in $ B_{r_{2}}(0),$ $
\eta\equiv 0$ in $ \mathbb{R}^{N}\setminus B_{r_{1}}(0)$ and
$|\nabla \eta|\leq 2/(r_{1}-r_{2}). $  As (\ref{3f}), we have
 \begin{eqnarray}\label{km}
&&\int_{B_{r_{1}}(0)}|x_{N}|^{2\alpha}\nabla u\nabla
(\eta^{2}\phi(u))\nonumber\\
 &\geq &C
\left(\int_{B_{r_{2}}(0)}|x_{N}|^{\alpha\cdot2^{*}(s)-s}|
h(u)|^{2^{*}(s)}\right)^{2/2^{*}(s)}
-\frac{4}{(r_{1}-r_{2})^{2}}\int_{B_{r_{1}}(0)}|x_{N}|^{2\alpha}h^{2}(u).
\end{eqnarray}
By (\ref{3c}) and  H\"older inequality,
\begin{eqnarray}\label{yh}
&&\int_{B_{r_{1}}(0)}|x_{N}|^{\alpha\cdot2^{*}(s)-s}u^{2^{*}(s)-1}\eta^{2}\phi(u)
\nonumber\\
&\leq&
{{t^{2}}\over{4(t-1)}}\int_{B_{r_{1}}(0)}|x_{N}|^{\alpha\cdot2^{*}(s)-s}|u|^{2^{*}(s)-2}|\eta
h(u)|^{2}\nonumber\\
&\leq& {{t^{2}}\over{4(t-1)}}\left(\int_{B_{r_{1}(0)}
}|x_{N}|^{\alpha\cdot2^{*}(s)-s}|u|^{2^{*}(s)t_{0}/2}\right)^{{2(2^{*}(s)-2)}
\over{2^{*}(s)t_{0}}}
\left(\int_{B_{r_{1}}(0)}|x_{N}|^{\alpha\cdot2^{*}(s)-s}| \eta h(
u)|^{2q}\right) ^{1/q},\nonumber\\
\end{eqnarray}
where $q={{ 2^{*}(s)t_{0}}\over{(t_{0}-2)2^{*}(s)+4}}$ satisfies
$q<2^{*}(s)/2,$ since $t_0>2$.
 Furthermore, by H\"older inequality, we have
\begin{eqnarray}\label{ggga}
\int_{B_{r_{1}}(0)}|x_{N}|^{2\alpha}h^{2}(u)\leq
\left(\int_{B_{r_{1}}(0)}
|x_{N}|^{(2\alpha-{{\alpha\cdot2^{*}(s)-s}\over{q}})q'}\right)^{1/q'}
\left(\int_{B_{r_{1}}(0)}|x_{N}|^{\alpha\cdot2^{*}(s)-s}|
h(u)|^{2q}\right)^{1/q},
\end{eqnarray}
where $ {1\over{q}}+{1\over{q'}}=1.$ By the fact that
$q\rightarrow 2^{*}(s)/2$ and $q'\rightarrow
2^{*}(s)/(2^{*}(s)-2)$
 as $t_{0}\rightarrow 2,$ we get that  $(2\alpha-{{\alpha\cdot2^{*}(s)-s}\over{q}})q'\rightarrow
2s/(2^{*}(s)-2)>0$ as $t_{0}\rightarrow 2.$ It follows that
$\int_{B_{1}(0)}
|x_{N}|^{(2\alpha-{{\alpha\cdot2^{*}(s)-s}\over{q}})q'}<\infty$ if
$t_{0}-2>0$ small enough. Then by  $(\ref{yy})- (\ref{ggga})$, we
have
\begin{eqnarray}\label{mjh}
\left(\int_{B_{r_{2}}}|x_{N}|^{\alpha\cdot2^{*}(s)-s}|h(u)|^{2^{*}(s)}\right)^{\frac{2}{2^{*}(s)}}
\leq
C'\left({{t^{2}}\over{4(t-1)}}+\frac{4}{(r_{1}-r_{2})^{2}}\right)
\left(\int_{B_{r_{1}}}|x_{N}|^{\alpha\cdot2^{*}(s)-s}| h(
u)|^{2q}\right) ^{\frac{1}{q}}.\nonumber
\end{eqnarray}
Letting $k\rightarrow\infty,$ we get
\begin{equation}\label{mjn}
|u|_{2^{*}(s)t/2,\ \alpha,s}\leq
C'^{\frac{1}{t}}\left({{t^{2}}\over{4(t-1)}}+\frac{4}{(r_{1}-r_{2})^{2}}\right)^{1/t}
|u|_{qt, \ \alpha,s},
\end{equation}
where $|u|_{l,\ \alpha,s}
:=(\int|x_{N}|^{\alpha\cdot2^{*}(s)-s}|u|^{l})^{1/l}.$ Choose
$\epsilon>0$ such that $(2+\epsilon)q<2^{*}(s)$. Let
$t_{n}=(2+\epsilon)(2^{*}(s)/2q)^{n-1},$
$r_{n}=\frac{\rho}{2}+(\frac{\rho}{2})^{n},$ $ n=1,2,\cdots.$ Then
by (\ref{mjn}) we have
$$|u|_{2^{*}(s)t_{n}/2,\ \alpha,s}\leq \left\{\prod^{n}_{i=2}
C'^{\frac{1}{t_{i}}}\left({{t^{2}_{i}}\over{4(t_{i}-1)}}+\frac{4}{(r_{i}-r_{i-1})^{2}}\right)^{1/t_{i}}\right\}
\cdot |u|_{(2+\epsilon)q, \ \alpha,s}.$$ Letting
$n\rightarrow\infty,$ we obtain that $u\in L^{\infty}(B_{\sigma})$
with $\sigma=\rho/2.$\hfill$\Box$

\medskip

\begin{Theorem}\label{nnn} (Harnack inequality)
Suppose that $\alpha>1/2$ and $0\leq s< 2.$ If $u\in
X_{\alpha}(B_{1}(0))$ is a  positive weak solution of equation
(\ref{yyyhsf2}), i.e., for every $\varphi\in
C^{\infty}_{0}(B_{1}(0)),$
$$\int_{B_{1}(0)}|x_{N}|^{2\alpha}\nabla u\nabla\varphi=
\int_{B_{1}(0)}|x_{N}|^{\alpha\cdot
2^{*}(s)-s}u^{2^{*}(s)-1}\varphi,$$ then there exist  constants
$C=C(N,s,\alpha)>0$ and $\varsigma=\varsigma(N,s,\alpha)\in (0,1)$
such that $$ \sup_{B_{\varsigma}(0)}u\leq
C\inf_{B_{\varsigma}(0)}u.$$

\end{Theorem}
\noindent{\bf Proof.} By the local boundedness of $u$ (Proposition
\ref{tererd}), $\forall\beta\in\mathbb{R}$ and $\eta\in
C^{\infty}_{0}(B_{1}(0))$, $\eta^{2}\overline{u}^{\beta}\in
X^{0}_{\alpha}(B_{1}(0)),$ where $\overline{u}=u+k$ and $k>0$.
 We have
\begin{eqnarray}\label{bc2}
\int_{B_{1}(0)}|x_{N}|^{2\alpha}\nabla
u\nabla(\eta^{2}\overline{u}^{\beta})=\beta
\int_{B_{1}(0)}|x_{N}|^{2\alpha}\eta^{2}\overline{u}^{\beta-1}|\nabla
\overline{u}|^{2}+2\int_{B_{1}(0)}|x_{N}|^{2\alpha}\eta
\overline{u}^{\beta}\nabla \overline{u}\nabla\eta.
\end{eqnarray}
Let's introduce a function $w$ defined by $w=\left\{
\begin{array}{l}
\overline{u}^{(\beta+1)/2}, \ \mbox{if}\  \beta\neq -1,\\
\log \overline{u},\quad \quad \mbox{if}\  \beta= -1.\\
\end{array} \right. $
Then we have
\begin{equation}\label{oiidffd}
\beta
\int_{B_{1}(0)}|x_{N}|^{2\alpha}\eta^{2}\overline{u}^{\beta-1}|\nabla
\overline{u}|^{2}=\left\{
\begin{array}{l}
\frac{4\beta}{(\beta+1)^{2}}\int_{B_{1}(0)}
|x_{N}|^{2\alpha}\eta^{2}|\nabla w|^{2}, \ \mbox{if}\  \beta\neq -1,\\
-\int_{B_{1}(0)}|x_{N}|^{2\alpha}\eta^{2}|\nabla w |^{2},\ \quad \quad \mbox{if}\  \beta= -1,\\
\end{array} \right.
\end{equation}
\begin{equation}\label{iiiqwws}
2\int_{B_{1}(0)}|x_{N}|^{2\alpha}\eta \overline{u}^{\beta}\nabla
\overline{u}\nabla\eta=\left\{
\begin{array}{l}
\frac{4}{\beta+1}\int_{B_{1}(0)}
|x_{N}|^{2\alpha}\eta w\nabla w\nabla \eta, \ \mbox{if}\  \beta\neq -1,\\
2\int_{B_{1}(0)}|x_{N}|^{2\alpha}\eta\nabla w \nabla\eta,\ \quad \quad \mbox{if}\  \beta= -1.\\
\end{array} \right.
\end{equation}
By $(\ref{bc2})-(\ref{iiiqwws})$, we obtain that if $\beta\neq
-1,0,$
\begin{eqnarray}\label{hhvcv}
&&\left|\int_{B_{1}(0)}|x_{N}|^{2\alpha}\nabla
\overline{u}\nabla(\eta^{2}\overline{u}^{\beta})\right|\nonumber\\
&=&\left|\frac{4\beta}{(\beta+1)^{2}}\int_{B_{1}(0)}|x_{N}|^{2\alpha}\eta^{2}|\nabla
w|^{2}+\frac{4}{\beta+1}\int_{B_{1}(0)}|x_{N}|^{2\alpha}\eta w\nabla w\nabla \eta\right|\nonumber\\
&\geq&
\frac{4|\beta|}{(\beta+1)^{2}}\int_{B_{1}(0)}|x_{N}|^{2\alpha}\eta^{2}|\nabla
w|^{2}\nonumber\\
&&-\frac{4}{|\beta+1|}\left(\frac{|\beta|}{2|\beta+1|}\int_{B_{1}(0)}|x_{N}|^{2\alpha}\eta^{2}|\nabla
w|^{2}+\frac{|\beta+1|}{2|\beta|}\int_{B_{1}(0)}|x_{N}|^{2\alpha}w^{2}|\nabla \eta|^{2}\right)\nonumber\\
&=&
\frac{2|\beta|}{(\beta+1)^{2}}\int_{B_{1}(0)}|x_{N}|^{2\alpha}\eta^{2}|\nabla
w|^{2}-\frac{2}{|\beta|}\int_{B_{1}(0)}|x_{N}|^{2\alpha}w^{2}|\nabla
\eta|^{2},
\end{eqnarray}
and if $\beta=-1,$
\begin{eqnarray}\label{hhvcvdd}
\left|\int_{B_{1}(0)}|x_{N}|^{2\alpha}\nabla
\overline{u}\nabla(\eta^{2}\overline{u}^{\beta})\right|
&=&\left|-\int_{B_{1}(0)}|x_{N}|^{2\alpha}\eta^{2}|\nabla w
|^{2}+2\int_{B_{1}(0)}|x_{N}|^{2\alpha}\eta\nabla w \nabla\eta\right|\nonumber\\
&\geq&\frac{1}{2}\int_{B_{1}(0)}|x_{N}|^{2\alpha}\eta^{2}|\nabla w
|^{2}-8\int_{B_{1}(0)}|x_{N}|^{2\alpha} |\nabla\eta|^{2}.
\end{eqnarray}
Moreover,  if $\beta\neq -1$, by $u\in
 L^{\infty}(B_{1}(0))$, we have
\begin{eqnarray}\label{ydfdf6}
\int_{B_{1}(0)}|x_{N}|^{\alpha\cdot
2^{*}(s)-s}u^{2^{*}(s)-1}\cdot\eta^{2}\overline{u}^{\beta}
&\leq&\int_{B_{1}(0)}|x_{N}|^{\alpha\cdot
2^{*}(s)-s}u^{2^{*}(s)-2}\cdot\eta^{2}w^{2}\nonumber\\
&\leq&C\int_{B_{1}(0)}|x_{N}|^{\alpha\cdot
2^{*}(s)-s}\cdot\eta^{2}w^{2}
\end{eqnarray}
and if $\beta=-1,$ by $\alpha\cdot 2^{*}(s)-s>\frac{1}{2}\cdot
2-s>-1$ and $u\in
 L^{\infty}(B_{1}(0))$, we obtain
\begin{eqnarray}\label{yf6}
\int_{B_{1}(0)}|x_{N}|^{\alpha\cdot
2^{*}(s)-s}u^{2^{*}(s)-1}\cdot\eta^{2}\overline{u}^{\beta}
&\leq&\int_{B_{1}(0)}|x_{N}|^{\alpha\cdot
2^{*}(s)-s}u^{2^{*}(s)-2} \nonumber\\ &\leq&
C\int_{B_{1}(0)}|x_{N}|^{\alpha\cdot 2^{*}(s)-s}<\infty.
\end{eqnarray}
By (\ref{hhvcvdd}), (\ref{yf6}) and
$\int_{B_{1}(0)}|x_{N}|^{2\alpha}\nabla
\overline{u}\nabla(\eta^{2}\overline{u}^{\beta})=\int_{B_{1}(0)}|x_{N}|^{\alpha\cdot
2^{*}(s)-s}u^{2^{*}(s)-1}\cdot\eta^{2}\overline{u}^{\beta}$, we
obtain that if $\beta=-1,$ then
\begin{equation}\label{cvcvc1}
\int_{B_{1}(0)}|x_{N}|^{2\alpha}\eta^{2}|\nabla w |^{2}\leq
16\int_{B_{1}(0)}|x_{N}|^{2\alpha}
|\nabla\eta|^{2}+2C\int_{B_{1}(0)}|x_{N}|^{\alpha\cdot
2^{*}(s)-s}<\infty.
\end{equation}
This means that for $\beta=-1,$ $w\in X_{\alpha, loc}(B_{1}(0)).$
Since
\begin{eqnarray}\label{dasa}
&&\int_{B_{1}(0)}|x_{N}|^{2\alpha}\eta^{2}|\nabla
w|^{2}\nonumber\\
&=& \int_{B_{1}(0)}|x_{N}|^{2\alpha}|\nabla (\eta w)|^{2}-2
\int_{B_{1}(0)}|x_{N}|^{2\alpha}w\nabla(\eta w)\nabla\eta-
\int_{B_{1}(0)}|x_{N}|^{2\alpha}w^{2}|\nabla\eta|^{2}\nonumber\\
&\geq&\frac{1}{2}\int_{B_{1}(0)}|x_{N}|^{2\alpha}|\nabla (\eta
w)|^{2}-3\int_{B_{1}(0)}|x_{N}|^{2\alpha}w^{2}|\nabla\eta|^{2},\nonumber
\end{eqnarray}
by $(\ref{hhvcv})$, $(\ref{ydfdf6})$ and
$\int_{B_{1}(0)}|x_{N}|^{2\alpha}\nabla
\overline{u}\nabla(\eta^{2}\overline{u}^{\beta})=\int_{B_{1}(0)}|x_{N}|^{\alpha\cdot
2^{*}(s)-s}u^{2^{*}(s)-1}\cdot\eta^{2}\overline{u}^{\beta}$, we
obtain that if $\beta\neq -1,0,$ then
\begin{eqnarray}\label{zzz9i9uhuhh}
&&\int_{B_{1}(0)}|x_{N}|^{2\alpha}|\nabla (\eta w)|^{2}
\nonumber\\
&\leq&
\frac{C(\beta+1)^{2}}{|\beta|}\int_{B_{1}(0)}|x_{N}|^{\alpha\cdot
2^{*}(s)-s}\eta^{2}w^{2}+\frac{(\beta+1)^{2}}{|\beta|^{2}}\int_{B_{1}(0)}|x_{N}|^{2\alpha}
w^{2}|\nabla\eta|^{2}.
\end{eqnarray}
When $2\alpha\geq\alpha\cdot 2^{*}(s)-s,$ by (\ref{zzz9i9uhuhh})
and Theorem \ref{yyyrtrtww}, we obtain that if $\beta\neq -1,0,$
\begin{eqnarray}\label{xxx12}
&&(\int_{B_{1}(0)}|x_{N}|^{\alpha\cdot 2^{*}(s)-s}|\eta
w|^{2^{*}(s)})^{2/2^{*}(s)}\nonumber\\&
\leq&\frac{C|\beta+1|^{2}}{|\beta|}\int_{B_{1}(0)}|x_{N}|^{\alpha\cdot
2^{*}(s)-s}\eta^{2}w^{2}+\frac{C|\beta+1|^{2}}{|\beta|^{2}}\int_{B_{1}(0)}|x_{N}|^{\alpha\cdot
2^{*}(s)-s} w^{2}|\nabla\eta|^{2}.\nonumber\\
\end{eqnarray}
When $2\alpha<\alpha\cdot 2^{*}(s)-s,$ we can choose
$s_{\alpha}\in (0,2)$ such that $2\alpha=\alpha\cdot
2^{*}(s_{\alpha})-s_{\alpha}.$ Then by (\ref{zzz9i9uhuhh}) and
Theorem \ref{yyyrtrtww}, we obtain that if $\beta\neq -1,0,$
\begin{eqnarray}\label{xxx12dc}
&&(\int_{B_{1}(0)}|x_{N}|^{2\alpha}|\eta
w|^{2^{*}(s_{\alpha})})^{2/2^{*}(s_{\alpha})}\nonumber\\
&=&(\int_{B_{1}(0)}|x_{N}|^{\alpha\cdot
2^{*}(s_{\alpha})-s_{\alpha}}|\eta
w|^{2^{*}(s_{\alpha})})^{2/2^{*}(s_{\alpha})}\leq
C\int_{B_{1}(0)}|x_{N}|^{2\alpha}|\nabla (\eta w)|^{2}\nonumber\\
&\leq&\frac{C|\beta+1|^{2}}{|\beta|}\int_{B_{1}(0)}|x_{N}|^{2\alpha}\eta^{2}w^{2}
+\frac{C|\beta|+1|^{2}}{|\beta|^{2}}\int_{B_{1}(0)}|x_{N}|^{2\alpha}
w^{2}|\nabla\eta|^{2}.
\end{eqnarray}
Let $r_{1}, r_{2}$ be such that $0<r_{1}<r_{2}<1$.  Let $\eta$ be
a cut-off function satisfying $\eta\equiv 1$ in $B_{r_{1}}(0)$,
$\eta\equiv 0$ in $\mathbb{R}^{N}\setminus B_{r_{2}}(0)$ and
$|\nabla\eta|\leq 2/(r_{2}-r_{1}).$ By (\ref{xxx12}) and
(\ref{xxx12dc}), we obtain that if $2\alpha\geq\alpha\cdot
2^{*}(s)-s,$
\begin{eqnarray}\label{xxx1asa2}
(\int_{B_{r_{1}}(0)}|x_{N}|^{\alpha\cdot 2^{*}(s)-s}|
w|^{2^{*}(s)})^{2/2^{*}(s)}\leq
\frac{C|1+\beta|}{(r_{2}-r_{1})^{2}}\int_{B_{r_{2}}(0)}|x_{N}|^{\alpha\cdot
2^{*}(s)-s}w^{2},
\end{eqnarray}
and if  $2\alpha<\alpha\cdot 2^{*}(s)-s,$
\begin{eqnarray}\label{xxaasa2}
(\int_{B_{r_{1}}(0)}|x_{N}|^{2\alpha}|
w|^{2^{*}(s_{\alpha})})^{2/2^{*}(s_{\alpha})}\leq
\frac{C|1+\beta|}{(r_{2}-r_{1})^{2}}\int_{B_{r_{2}}(0)}|x_{N}|^{2\alpha}w^{2},
\end{eqnarray}
where $C>0$ is a constant depending only on $\beta$ and is bounded
when $|\beta|$ is bounded away from zero.

Set $\gamma=1+\beta$ and set $\Phi(p,r)=\left\{
\begin{array}{l}
(\int_{B_{r}(0)}|x_{N}|^{\alpha\cdot 2^{*}(s)-s}|
\overline{u}|^{p})^{1/p}, \ \mbox{if}\ 2\alpha\geq\alpha\cdot 2^{*}(s)-s\\
(\int_{B_{r}(0)}|x_{N}|^{2\alpha}|
\overline{u}|^{p})^{1/p}, \quad\quad\quad \mbox{if}\ 2\alpha<\alpha\cdot 2^{*}(s)-s.\\
\end{array} \right.$
By (\ref{xxx1asa2}) and (\ref{xxaasa2}), we obtain that when
$2\alpha\geq\alpha\cdot 2^{*}(s)-s,$
\begin{equation}\label{yerde1}
\Phi\left(\frac{2^{*}(s)}{2}\gamma,r_{1}\right)\leq
\left(\frac{C(1+|\gamma|)}{(r_{2}-r_{1})^{2}}\right)^{1/|\gamma|}
\Phi(\gamma,r_{2}),\ \mbox{if}\ \gamma>0,
\end{equation}
\begin{equation}\label{yzze1}
\Phi(\gamma,r_{2})\leq
\left(\frac{C(1+|\gamma|)}{(r_{2}-r_{1})^{2}}\right)^{1/|\gamma|}
\Phi\left(\frac{2^{*}(s)}{2}\gamma,r_{1}\right),\ \mbox{if}\
\gamma<0.
\end{equation}
and when $2\alpha<\alpha\cdot 2^{*}(s)-s,$
\begin{equation}\label{yerde1s}
\Phi\left(\frac{2^{*}(s_{\alpha})}{2}\gamma,r_{1}\right)\leq
\left(\frac{C(1+|\gamma|)}{(r_{2}-r_{1})^{2}}\right)^{1/|\gamma|}
\Phi(\gamma,r_{2}),\ \mbox{if}\ \gamma>0,
\end{equation}
\begin{equation}\label{yzze1fd}
\Phi(\gamma,r_{2})\leq
\left(\frac{C(1+|\gamma|)}{(r_{2}-r_{1})^{2}}\right)^{1/|\gamma|}
\Phi\left(\frac{2^{*}(s_{\alpha})}{2}\gamma,r_{1}\right),\
\mbox{if}\ \gamma<0,
\end{equation}
 Hence taking $
p>0$, we set $\gamma=\gamma_{m}=p(2^{*}(s)/2)^{m-1},$ and for
$\varsigma\in (0,1/4)$, set
$r_{m}=\varsigma+(\frac{\varsigma}{4})^{m},$ $m=1,2,\cdots,$ so
that by (\ref{yerde1}) or (\ref{yerde1s}), $
\Phi(\gamma_{m},\varsigma) \leq C\Phi(p,5\varsigma/4),$ $
m=1,2,\cdots$.
 Consequently, letting $m$ tend
to infinity, we have
\begin{equation}\label{idgdf}
\sup_{B_{\varsigma}(0)}\overline{u}\leq
C\Phi\left(p,\frac{5\varsigma}{4}\right).
\end{equation}
In a similar manner, by (\ref{yzze1}) or (\ref{yzze1fd}), we can
prove that for any $p>0,$
\begin{equation}\label{xwqqw}
\Phi\left(-p,\frac{5\varsigma}{4}\right)\leq
C\Phi(-\infty,\varsigma)=C\inf_{B_{\varsigma}(0)}\overline{u}.
\end{equation}
Let $S_{\rho}(x,r)$ be the ball $\{y\in\mathbb{R}^N\ |\
\rho(x,y)<r\}$ with the metric $\rho$ defined in \cite[Theorem
2.7]{FL}. By Proposition 2.9 and Theorem 2.7 of \cite{FL}, we
deduce that
  when $\varsigma$ small enough, there exists
 $\delta>0$ such that
 \begin{eqnarray}\label{vvcbc77d6dqq}
 B_{5\varsigma/4}(0)\subset S_{\rho}(0,\delta)\subset
 B_{1/2}(0).\end{eqnarray}
By (\ref{cvcvc1}), \cite[Theorem 3.2 and Remark 3.3]{FL}(see page
538 and 539) and (4.2f) of \cite{FL}(see page 538), we deduce that
there exist $p>0$ and constant $C>0$ such that
$$(\int_{S_{\rho}(0,\delta)}\bar{u}^p)(\int_{S_{\rho}(0,\delta)}\bar{u}^{-p})\leq C.$$
Then by (\ref{vvcbc77d6dqq}), we get that
\begin{equation}\label{oqzaa}
\Phi(p,5\varsigma/4)\Big/\Phi(-p,5\varsigma/4) \leq
C.\end{equation} Letting $k\rightarrow 0$, by
$(\ref{idgdf})-(\ref{oqzaa})$, we get that
$\sup_{B_{\varsigma}(0)}u\leq
C\inf_{B_{\varsigma}(0)}u.$\hfill$\Box$

\medskip

Using the similar argument as the proofs of
 the above two theorems  and Theorem 8.22 of \cite{GL}, we can get the
 following theorem

\begin{Theorem}\label{czxzxz}
Suppose that $\alpha>1/2$ and $0\leq s< 2.$ If $u\in
X_{\alpha}(B_{1}(0))$ is a  weak solution of equation
(\ref{yyyhsf2}) in $B_{1}(0),$  then there exists $\sigma\in
(0,1)$ such that $u\in C^{0,\gamma}(\overline{B_{\sigma}(0)})$ for
some $0<\gamma<1.$
\end{Theorem}

\begin{proposition}\label{tqaaa3w}
If  the same conditions as Theorem \ref{czxzxz} holds, then
 $u\in C^{0,\gamma}(B_{1}(0))\cap
C^{2,\gamma}(B_{1}(0)\setminus \{x_{N}=0\})$ and $\frac{\partial
u}{\partial x_{i}}\in C^{0,\gamma}(B_{1}(0)),$ $1\leq i\leq N-1$
for some $0<\gamma<1.$
\end{proposition}
\noindent{\bf Proof.} By Theorem \ref{czxzxz} and Schauder
estimates, we obtain that $u\in C^{2,\gamma}(B_{1}(0)\setminus
\{x_{N}=0\}),$ since the operator $-div(|x_{N}|^{2\alpha}\nabla
u)$ is uniformly elliptic in compact subset of $B_{1}(0)\setminus
\{x_{N}=0\}.$ As the same proof of Theorem 8.8 in \cite{GL}, we
know that for $1\leq i\leq N-1,$ $\frac{\partial u}{\partial
x_{i}}\in X_{\alpha,loc}(B_{1}(0))$ and it satisfies
$-div\left(|x_{N}|^{2\alpha}\nabla \left(\frac{\partial
u}{\partial x_{i}}\right)\right) =(2^{*}(s)-1)|x_{N}|^{\alpha\cdot
2^{*}(s)-s}u^{2^{*}(s)-2}\frac{\partial u}{\partial x_{i}}$ in
$B_{1}(0)$ weakly. Using the same method as the proof of Theorem
\ref{czxzxz}, we can get  $\frac{\partial u}{\partial x_{i}}\in
C^{0,\gamma}(B_{1}(0)),$ $1\leq i\leq N-1$ for some
$0<\gamma<1.$\hfill$\Box$

\begin{proposition}\label{fff1fc}
Let $\Omega=\{(x',x_{N})\ |\  |x_{N}|\leq 2,\ |x'|\leq 1 \}$
 and $\Omega_{1}=\{(x',x_{N})\ |\  |x_{N}|\leq
\frac{1}{2},\ |x'|\leq \frac{1}{2} \}.$ If $\alpha>1/2$, $0\leq s<
2$ and $u$ is a weak solution of equation (\ref{yyyhsf2}) in
$\Omega$, then there exists  $C>0$ such that
$$\left|\frac{\partial u(x) }{\partial
x_{N}}\right|\leq C|x_{N}|^{-1},\ \forall x\in
\Omega_{1}\setminus\{x \ |\ x_{N}=0\}.$$
\end{proposition}
\noindent{\bf Proof.}
  For $0<\epsilon\leq1,$ let $\Omega_{\epsilon}=\{
(x',x_{N})\ |\ \frac{1}{4}\epsilon\leq x_{N}\leq 2\epsilon,\
|x'|\leq \epsilon\}$ and  $\Omega^{*}=\{(y',y_{N})\ |\
\frac{1}{4}\leq y_{N}\leq 2,\ |y'|\leq 1 \}.$ For $x\in
\Omega_{\epsilon},$ set $u_{\epsilon}(y)=u(\epsilon y),$
 $y=x/\epsilon.$ By Theorem \ref{tererd}, $u_{\epsilon}$ is bounded in $\Omega^{*}$. Straightforward
calculation shows that $u_{\epsilon}$ satisfies
$-div(|y_{N}|^{2\alpha}\nabla u_{\epsilon})
=\epsilon^{\beta}|u_{\epsilon}|^{2^{*}(s)-2}u_{\epsilon}$ in $
\Omega^{*}$, where $\beta=\alpha\cdot (2^{*}(s)-2)+2-s\geq0.$
$L^{p}-$ estimate gives that there exists $C'>0$ such that
\begin{eqnarray}\label{hhh231}
||u_{\epsilon}||_{C^{1,\gamma}(\Omega_{2})}\leq
C'(||u_{\epsilon}||_{L^{\infty}(\Omega^{*})}+\epsilon^{\beta}
||u^{2^{*}(s)-1}_{\epsilon}||_{L^{\infty}(\Omega^{*})}) \leq
C'(M+\epsilon^{\beta}M^{2^{*}(s)-1}) :=C,\nonumber
\end{eqnarray}
where $\Omega_{2}=\{(y',y_{N})\ |\ \frac{1}{2}< y_{N}\leq 1,\
|y'|\leq \frac{1}{2} \}$.
 In
particular, we have $\left|\epsilon\frac{\partial u}{\partial
x_{N}}(0,\epsilon)\right|=\left|\frac{\partial
u_{\epsilon}}{\partial y_{N}}(0,1)\right|\leq C.$ Thus
\begin{equation}\label{oq3}
\left|\frac{\partial u}{\partial x_{N}}(0,x_{N})\right|\leq
C|x_{N}|^{-1}, \ 0<x_{N}\leq 1.
\end{equation}
For fixed $|x'_{0}|\leq\frac{1}{2}$, consider
$\widetilde{u}(x)=u(x+(x'_{0},0))$. As (\ref{oq3}), we have $
\left|\frac{\partial u}{\partial x_{N}}(x'_{0},x_{N})\right|\leq
C|x_{N}|^{-1}.$ Therefore $\left|\frac{\partial u(x) }{\partial
x_{N}}\right|\leq C|x_{N}|^{-1},\ \forall x\in
\Omega_{1}\setminus\{x \ |\ x_{N}=0\}.$\hfill$\Box$

\section{Symmetry and uniqueness of solutions}\label{rft6678891}
\hspace*{\parindent} In this section, we obtain some symmetry and
uniqueness results for positive solutions of equations
(\ref{yyyhsf2}).

Given $\lambda>0$ and a function
$u:\mathbb{R}^{N}\rightarrow\mathbb{R}$, define $u_{\lambda}(x)
=\frac{\lambda^{N-2+2\alpha}}{|x|^{N-2+2\alpha}}u\left(\frac{\lambda^{2}x}{|x|^{2}}\right),\
x\in\mathbb{R}^{N}\setminus\{0\}.$ We shall use the method of
moving sphere (see \cite{CaoLi, LZ, LZh}) and its variant (see
\cite{MM}) to prove the following Theorem
\begin{Theorem}\label{yqwe12}  Suppose that $\alpha>1/2$ and $0\leq s<2$. If
$u\in X_{\alpha,loc}(\mathbb{R}^{N})$ is a positive solution of
equation (\ref{yyyhsf2}), then there exists a positive number
$\lambda$ such that $u(x)=u_{\lambda}(x),$
$x\in\mathbb{R}^{N}\setminus \{0\}.$
\end{Theorem}
\noindent{\bf Proof.} By Theorem \ref{czxzxz}, we know that $u$ is
H\"older continuous in $\mathbb{R}^{N}$. The proof is divided into
four steps.

\medskip

\underline{Step 1}. In this step, we shall prove that there exists
$\lambda_{0}>0$ such that $u_{\lambda}(x)\geq u(x),$
$|x|\leq\lambda$ if $0<\lambda<\lambda_{0}.$

By Proposition \ref{yhadsd}, we know that $u_{\lambda}$ satisfies
equation (\ref{yyyhsf2}) and  $u_{\lambda}|_{\partial
B_{\lambda}(0)}\equiv u|_{\partial B_{\lambda}(0)}.$ Thus
\begin{eqnarray}\label{uqasdx}
-div(|x_{N}|^{2\alpha}\nabla(u_{\lambda}-u))
&=&|x_{N}|^{\alpha\cdot
2^{*}(s)-s}(u^{2^{*}(s)-1}_{\lambda}-u^{2^{*}(s)-1})\nonumber\\
&=&(2^{*}(s)-1)|x_{N}|^{\alpha\cdot 2^{*}(s)-s}\psi_{\lambda}^{
2^{*}(s)-2}(u_{\lambda}-u),
\end{eqnarray}
where $\psi_{\lambda}(x)$ is some number between $u_{\lambda}(x)$
and $u(x)$.  Let $\Omega^{-}_{\lambda}=\{x\in B_{\lambda}(0)\ |\
(u_{\lambda}-u)(x)\leq 0\}.$ Set $M=\max_{B_{1}(0)}u.$ By Theorem
\ref{tererd}, we have $M<+\infty.$ Multiplying equation
(\ref{uqasdx}) by $(u_{\lambda}-u)^{-}$ and integrating. By
H\"older inequality and Theorem \ref{yyyrtrtww}, we get that if
$0<\lambda\leq 1,$ then
\begin{eqnarray}\label{yq12}
&&\int_{\Omega^{-}_{\lambda}}|x_{N}|^{2\alpha}|\nabla(u_{\lambda}-u)|^{2}\nonumber\\
& =&(2^{*}(s)-1)\int_{\Omega^{-}_{\lambda}}|x_{N}|^{\alpha\cdot
2^{*}(s)-s}\psi_{\lambda}^{
2^{*}(s)-2}(u_{\lambda}-u)^{2}\nonumber\\
&\leq&
(2^{*}(s)-1)\left(\int_{\Omega^{-}_{\lambda}}|x_{N}|^{\alpha\cdot
2^{*}(s)-s}\psi_{\lambda}^{
2^{*}(s)}\right)^{\frac{2^{*}(s)-2}{2^{*}(s)}}
\left(\int_{\Omega^{-}_{\lambda}}|x_{N}|^{\alpha\cdot
2^{*}(s)-s}(u_{\lambda}-u)^{2^{*}(s)}\right)^{\frac{2}{2^{*}(s)}}\nonumber\\
&\leq&(2^{*}(s)-1)M^{2^{*}(s)-2}\left(\int_{\Omega^{-}_{\lambda}}|x_{N}|^{\alpha\cdot
2^{*}(s)-s}\right)^{\frac{2^{*}(s)-2}{2^{*}(s)}}
\left(\int_{\Omega^{-}_{\lambda}}|x_{N}|^{\alpha\cdot
2^{*}(s)-s}(u_{\lambda}-u)^{2^{*}(s)}\right)^{\frac{2}{2^{*}(s)}}\nonumber\\
&\leq&
CM^{2^{*}(s)-2}\left(\int_{\Omega^{-}_{\lambda}}|x_{N}|^{\alpha\cdot
2^{*}(s)-s}\right)^{\frac{2^{*}(s)-2}{2^{*}(s)}}
\int_{\Omega^{-}_{\lambda}}|x_{N}|^{2\alpha}|\nabla(u_{\lambda}-u)|^{2}.
\end{eqnarray}
Since $\lim_{\lambda\rightarrow
0}\int_{\Omega^{-}_{\lambda}}|x_{N}|^{\alpha\cdot 2^{*}(s)-s}=0,$
by (\ref{yq12}), we deduce that
$\int_{\Omega^{-}_{\lambda}}|x_{N}|^{2\alpha}|\nabla(u_{\lambda}-u)|^{2}=0$
if $\lambda>0$ small enough. It follows that if $\lambda$ small
enough, then  for any $x\in B_{\lambda}(0)$, $u_{\lambda}(x)\geq
u(x)$.

\medskip

\underline{Step 2}. Set $\overline{\lambda}=\sup\{\mu>0\ |\
u_{\lambda}(x)\geq u(x), \ |x|\leq\lambda,\ 0<\lambda<\mu \}.$ We
shall prove that if $\overline{\lambda}<\infty,$ then
$u_{\overline{\lambda}}\equiv u$ in $\mathbb{R}^{N}\setminus
\{0\}.$

\medskip

Obviously, it is sufficient to prove that
$u_{\overline{\lambda}}\equiv u$ in $B_{\overline{\lambda}}(0).$
From the definition of $\overline{\lambda}$, we know that
$u_{\overline{\lambda}}\geq u$ in $B_{\overline{\lambda}}(0).$ If
$u_{\overline{\lambda}}\not\equiv u$ in
$B_{\overline{\lambda}}(0),$ by
\begin{eqnarray}\label{tqzx}
-div(|x_{N}|^{2\alpha}\nabla(u_{\overline{\lambda}}-u))
=(2^{*}(s)-1)|x_{N}|^{\alpha\cdot
2^{*}(s)-s}\psi_{\overline{\lambda}}^{
2^{*}(s)-2}(u_{\overline{\lambda}}-u)\geq 0,\nonumber
\end{eqnarray}
and Proposition \ref{mvxcee}, we get that
\begin{equation}\label{zmd}
(u_{\overline{\lambda}}-u)(x)>0,\ \forall x\in
B_{\overline{\lambda}}(0).
\end{equation}
It follows that for  $\delta>0$ small enough,
\begin{eqnarray}\label{azsuqads}
\max_{x\in\partial
B_{\overline{\lambda}-\delta}(0)}(u_{\overline{\lambda}}-u)(x)>0.
\end{eqnarray}
By (\ref{azsuqads}), we can choose $\epsilon=\epsilon(\delta)>0$
small enough, such that $\epsilon=o(\delta)$ as $\delta\rightarrow
0$ and
\begin{eqnarray}\label{uqads}
\max_{x\in\partial B_{\lambda-\delta}(0)}(u_{\lambda}-u)(x)>0,\
\mbox{if}\
\overline{\lambda}\leq\lambda\leq\overline{\lambda}+\epsilon.
\end{eqnarray}

Set $\Omega^{-}_{\lambda}=\{x\ |\ \lambda-\delta\leq
|x|\leq\lambda,\ u_{\lambda}(x)-u(x)\leq0\}.$ By (\ref{uqads}) and
the fact that $(u_{\lambda}-u)|_{\partial B_{\lambda}(0)}\equiv
0$, we get that $(u_{\lambda}-u)^{-}|_{\partial
\Omega^{-}_{\lambda}}\equiv0.$ Then as (\ref{yq12}), we get
\begin{eqnarray}\label{yq12v}
&&\int_{\Omega^{-}_{\lambda}}|x_{N}|^{2\alpha}|\nabla(u_{\lambda}-u)|^{2}\nonumber\\
& =&(2^{*}(s)-1)\int_{\Omega^{-}_{\lambda}}|x_{N}|^{\alpha\cdot
2^{*}(s)-s}\psi_{\lambda}^{
2^{*}(s)-2}(u_{\lambda}-u)^{2}\nonumber\\
&\leq&
CM^{2^{*}(s)-2}\left(\int_{\Omega^{-}_{\lambda}}|x_{N}|^{\alpha\cdot
2^{*}(s)-s}\right)^{\frac{2^{*}(s)-2}{2^{*}(s)}}
\int_{\Omega^{-}_{\lambda}}|x_{N}|^{2\alpha}|\nabla(u_{\lambda}-u)|^{2}.
\end{eqnarray}
By
$$\lim_{\delta\rightarrow 0}\int_{\Omega^{-}_{\lambda}}|x_{N}|^{\alpha\cdot
2^{*}(s)-s}=0$$ and (\ref{yq12v}), we know that if $\delta$ small
enough,  Lebesgue measure of $\Omega^{-}_{\lambda}$ must be zero.
Thus when
$\overline{\lambda}\leq\lambda\leq\overline{\lambda}+\epsilon,$
\begin{equation}\label{bzcxx}
u_{\lambda}(x)-u(x)\geq 0,\ \lambda-\delta\leq |x|\leq\lambda.
\end{equation}

By (\ref{zmd}), we deduce  that there exists $C(\delta)>0$ such
that $u_{\overline{\lambda}}(x)-u(x)\geq C(\delta)>0$ if
$|x|<\lambda-\delta$.  Thus we can choose $\epsilon$ small enough,
such that $u_{\lambda}(x)-u(x)>0$ if $|x|<\lambda-\delta$ and
$\overline{\lambda}\leq\lambda\leq\overline{\lambda}+\epsilon.$
Combining (\ref{bzcxx}), we obtain that $u_{\lambda}(x)-u(x)\geq0$
 if  $|x|\leq\lambda$ and
$\overline{\lambda}\leq\lambda\leq\overline{\lambda}+\epsilon$. It
contradicts the definition of $\overline{\lambda}.$ Thus
$u_{\overline{\lambda}}\equiv u$ in $\mathbb{R}^{N}\setminus
\{0\}.$

\medskip

\underline{Step 3}. For $b\in\mathbb{R}^{N-1},$ let
$u^{(b)}(x)=u(x+(b,0)),$ $x\in \mathbb{R}^{N}$ and let
$\overline{\lambda}_{b}$ be defined as in Step 2 relative to
$u^{(b)}$. In this step, we shall prove that if
$\overline{\lambda}_{b}=\infty$ for some $b\in\mathbb{R}^{N-1},$
then $\overline{\lambda}_{b}=\infty$ for all
$b\in\mathbb{R}^{N-1}.$

\medskip

By Step 2, there is a maximal $\overline{\lambda}_{b}>0$ such that
$(u^{(b)})_{\lambda}(x)\geq u^{(b)}(x),$ if  $ |x|\leq\lambda$ and
$0<\lambda<\overline{\lambda}_{b}$. It follows that
$(u^{(b)})_{\lambda}(x)\leq u^{(b)}(x),$ if $ |x|\geq\lambda$ and
$0<\lambda<\overline{\lambda}_{b}$. Letting $x_{b}=x-(b,0),$ we
have $u(x)\geq \left(\frac{\lambda}{|x_{b}|}\right)^{N-2+2\alpha}
u\left(\frac{\lambda^{2}x_{b}}{|x_{b}|^{2}}+(b,0)\right).$ Since
$\overline{\lambda}_{b}=\infty,$ we know that the above inequality
holds for all $\lambda>0$ and $|x_{b}|\geq\lambda.$ For any fixed
$\lambda>0,$ it follows that
$$
\lim_{|x|\rightarrow\infty}|x|^{N-2+2\alpha}u(x)
\geq\lim_{|x_{b}|\rightarrow\infty}\left(\frac{\lambda|x|}{|x_{b}|}\right)^{N-2+2\alpha}
u\left(\frac{\lambda^{2}x_{b}}{|x_{b}|^{2}}+(b,0)\right)
=\lambda^{N-2+2\alpha}u(b,0).
$$
Letting $\lambda\rightarrow\infty,$ this implies
$\lim_{|x|\rightarrow\infty}|x|^{N-2+2\alpha}u(x)=\infty.$ Assume
that there is $a\neq b$ such that $\overline{\lambda}_{a}<\infty.$
Then by Step 2, it is $u^{(a)}\equiv
(u^{(a)})_{\overline{\lambda}_{a}},$ i.e., $u(x)=
\left(\frac{\lambda}{|x_{a}|}\right)^{N-2+2\alpha}
u\left(\frac{\lambda^{2}x_{a}}{|x_{a}|^{2}}+(a,0)\right),$
$x_{a}=x-(a,0).$ This gives
$\displaystyle\lim_{|x|\rightarrow\infty}|x|^{N-2+2\alpha}u(x)=\overline{\lambda}_{a}u(a,0)<\infty,$
which contradicts to
$\displaystyle\lim_{|x|\rightarrow\infty}|x|^{N-2+2\alpha}u(x)=\infty.$

\medskip

\underline{Step 4}. In this step, we shall prove that
$\overline{\lambda}_{b}<\infty$ for all $ b\in\mathbb{R}^{N-1}.$

\medskip

By contradiction. If not, then by Step 3, for any
$b\in\mathbb{R}^{N-1},$ $\overline{\lambda}_{b}=\infty.$ Let
$$g_{\lambda,b}(x)=u(x+(b,0))-\left(\frac{\lambda}{|x|}\right)^{N-2+2\alpha}
u\left(\frac{\lambda^{2}x}{|x|^{2}}+(b,0)\right)$$ for
$|x|\leq\lambda,$ $\lambda>0$ and $b\in\mathbb{R}^{N-1}.$ Then $
g_{|x|,b}(rx)=u(rx+(b,0))-\frac{1}{r^{N-2+2\alpha}}u\left(\frac{x}{r}+(b,0)\right).
$ Since $g_{\lambda,b}(x)<0$ if $|x|\leq\lambda$,
 we get that when $0<r<1,$  \begin{equation}\label{cxz}
 g_{|x|,b}(rx)<0.\end{equation}
By Proposition \ref{tqaaa3w}, we deduce that for any $x\neq 0,$
$\frac{d}{dr}(g_{|x|,b}(rx))|_{r=1}$ exists. By straightforward
calculation, we have
\begin{eqnarray}\label{xzsaw}
&&\frac{d}{dr}(g_{|x|,b}(rx))\Big|_{r=1}\nonumber\\
 &=&
\left(x(\nabla u)(rx+(b,0))+\frac{x}{r^{N+2\alpha}}(\nabla
u)(\frac{x}{r}+(b,0))
+\frac{N-2+2\alpha}{r^{N-1+2\alpha}}u(\frac{x}{r}+(b,0))\right)\Big|_{r=1}\nonumber\\
&=& 2x\nabla u(x+(b,0))+(N-2+2\alpha)u(x+(b,0)).
\end{eqnarray}
By (\ref{cxz}) and the fact that $g_{|x|,b}(x)=0,$ we get that $
\frac{d}{dr}(g_{|x|,b}(rx))|_{r=1}\geq 0. $ Thus by (\ref{xzsaw}),
we have
\begin{equation}\label{cxxqqqw}
2(x-(b,0))\nabla u(x)+(N-2+2\alpha)u(x)\geq 0
\end{equation}
Divided both side by $|b|$ in (\ref{cxxqqqw}) and let $|b|$ tend
to $\infty,$ we get that $a\nabla u_{x'}(x)\leq 0$ for any
$a\in\mathbb{R}^{N-1}$ with  $|a|=1.$
 It follows that $ \nabla u_{x'}(x)=0,\
\forall x\in\mathbb{R}^{N}. $ Thus $u$ is independent of $x',$
i.e., $u=u(x_{N}).$ By Theorem \ref{tererd}, we get that
$u=u(x_{N})$ is a positive solution of
$-(x^{2\alpha}_{N}u')'=x^{\alpha\cdot
2^{*}(s)-s}_{N}u^{2^{*}(s)-1}$, $x_{N}>0$ with
$u(0)=\lim_{x_{N}\rightarrow 0}u(x_{N})<\infty.$ However, by Lemma
\ref{yyrhfbvgg}, we know that this equation has no positive
solution. Thus $\overline{\lambda}_{b}<\infty$ for any $b\in
\mathbb{R}^{N-1}$. \hfill$\Box$

\begin{remark}\label{ccv5}
By this Theorem,
 we know that if  $u\in X_{\alpha,loc}(\mathbb{R}^{N})$
 is a positive solution of equation (\ref{yyyhsf2}),
then $u(x)=O(|x|^{-(N-2+2\alpha)}).$ It follows that  if $u\in
X_{\alpha,loc}(\mathbb{R}^{N})$ is a positive solution of equation
(\ref{yyyhsf2}), then $u\in X_{\alpha}(\mathbb{R}^{N})$.
 \end{remark}

 \begin{lemma}\label{yyrhfbvgg} If $\alpha>1/2$ and $0\leq s<2,$
 then the following equation \begin{equation}\label{uu7tgrf44bbvg}
 -(r^{2\alpha}f')'=r^{\alpha\cdot 2^{*}(s)-s}f^{2^{*}(s)-1},\
 r>0,
 \ f(r)>0,\ f(0)=\lim_{r\rightarrow 0+}f(r)<\infty
 \end{equation} has no
 solution.
 \end{lemma}
 \noindent{\bf Proof.} Equation (\ref{uu7tgrf44bbvg}) is equivalent to
$ -f''(r)-\frac{2\alpha}{r}f'(r)=r^{\beta}f^{2^{*}(s)-1}(r),\ r>0,
$ where $\beta=\alpha\cdot 2^{*}(s)-s-2\alpha.$ For $\tau>0,$
making the change of variable $y=r^{\tau}/\tau,$ $f(r)=u(y)$, this
equation  becomes
\begin{equation}\label{9oowtdffvc}
-u''(y)-\frac{\frac{2\alpha-1}{\tau}+1}{y}u(y)=
\tau^{\sigma}y^{\sigma}u^{2^{*}(s)-1}(y),\ y>0,
\end{equation}
where $\sigma=\frac{\beta-2(\tau-1)}{\tau}$.
 For $x\in\mathbb{R},$ let $[x]$ denote the
 largest one among the integers which do not exceed $x$. Let
 $\tau=(2\alpha-1)/[2\alpha]$ and
 $k=\frac{2\alpha-1}{\tau}+2=[2\alpha]+2$. Then equation
 (\ref{9oowtdffvc}) is equivalent to
 \begin{equation}\label{99iqnnanzggv}
 -\triangle u(x)=\tau^{\sigma}|x|^{\sigma}
 u^{2^{*}(s)-1}(x),\ x\in \mathbb{R}^{k},
 \end{equation}
where  $|x|=y.$ It is easy to verify that $\sigma>-2$ and
$\frac{k+2+2\sigma}{k-2}>2^{*}(s)-1.$ Thus according to
Proposition 5.2 in \cite{SZ} by Serrin and Zou, radially symmetric
positive solutions of (\ref{99iqnnanzggv}), if any, satisfy
$\lim_{|x|\rightarrow 0}
|x|^{\frac{2+\sigma}{2^{*}(s)-2}}u(x)=\lambda,$ for some positive
constant $\lambda.$ This contradicts the assumption $
f(0)=\lim_{r\rightarrow 0+}f(r)<\infty.$ This lemma is
established. \hfill$\Box$

\medskip

By  Remark \ref{ccv5}, Proposition \ref{jjjhgrt}
 and the classical moving plane method (see \cite{GNN}), we can get the following
 Theorem
 \begin{Theorem}\label{vxcsde} Suppose that $\alpha>1/2$ and $0\leq s<2.$
 If $u\in X_{\alpha,loc}(\mathbb{R}^{N})$ is a positive
solution of equation (\ref{yyyhsf2}), then there exists
$x'_{0}\in\mathbb{R}^{N-1}$ such that $u$ is axially symmetric
about the axis $\{x=(x',x_{N})\in\mathbb{R}^{N}\ |\ x'=x'_{0}\},$
i.e., $u(x',x_{N})=u(|x'-x'_{0}|,x_{N})$. Moreover,
$\frac{\partial u}{\partial r}(r,x_{N})<0$ for $r=|x'-x'_{0}|>0.$
 \end{Theorem}

\begin{Theorem}\label{uuxxxc}
If $u\in X_{\alpha,loc}(\mathbb{R}^{N})$ is a nonnegative weak
solution of $div(|x_{N}|^{2\alpha}\nabla u)=0$  in $
\mathbb{R}^{N}$, then $u\equiv a$ for some constant $a\geq 0.$
\end{Theorem}
\noindent{\bf Proof.} As the proof of Theorem \ref{yqwe12}, we can
get that either $u\equiv u_{\overline{\lambda}}$ for some
$0<\overline{\lambda}<\infty,$
  or $\overline{\lambda}_{b}=\infty$ for
all $b\in\mathbb{R}^{N-1}$,  where $\overline{\lambda}_{b}$ is
defined in the Step 3 of proof Theorem \ref{yqwe12}. If
$\overline{\lambda}_{b}=\infty$ for all $b\in\mathbb{R}^{N}$, then
as  step 4 of the proof of Theorem \ref{yqwe12}, we deduce  that
$u\equiv a$ for some constant $a\geq 0.$ If $u\equiv
u_{\overline{\lambda}}$ for some $0<\overline{\lambda}<\infty,$
then $u(x)=O(|x|^{-(N-2+2\alpha)})$ as $|x|\rightarrow\infty.$
Then by Proposition \ref{rwewe}, we get that  $ 0\leq
\sup_{B_{R}(0)}u\leq \sup_{\partial B_{R}(0)}u\rightarrow 0,$ as
$R\rightarrow\infty.$ It follows that $u\equiv 0.$\hfill$\Box$

\begin{Theorem}\label{bz2wdf} Suppose that $\alpha>1/2$ and $0\leq s<2.$
Let $u\in X_{\alpha,loc}(\mathbb{R}^{N})$ be a positive solution
of equation (\ref{yyyhsf2}) which satisfies
$$u(x)=|x|^{-(N-2+2\alpha)}u(x/|x|^{2}):=\widetilde{u}(x),\ \forall x\in
\mathbb{R}^{N}\setminus\{0\}.$$ Then there is
$x'_{0}\in\mathbb{R}^{N-1}$ such that
$u(x',0)=u(x'_{0},0)(1+|x'-x'_{0}|^{2})^{-(N-2+2\alpha)/2},$
$\forall x'\in\mathbb{R}^{N-1}.$
\end{Theorem}
\noindent{\bf Proof.}  For a fixed $b\in\mathbb{R}^{N-1},$ define
$u^{(b)}(x)=u(x+(b,0)).$ By Theorem \ref{yqwe12}, there exists
$\lambda_{b}>0$ such that $u^{(b)}=(u^{(b)})_{\lambda_{b}},$ i.e.,
$u(x+(b,0))=\left(\frac{\lambda_{b}}{|x|}\right)^{N-2+2\alpha}
u\left(\frac{\lambda^{2}_{b}x}{|x|^{2}}+(b,0)\right).$ Letting
$x_{b}=x-(b,0)$ for all $x,$ this identity becomes
\begin{equation}\label{xcretyu}
u(x)=\left(\frac{\lambda_{b}}{|x_{b}|}\right)^{N-2+2\alpha}
u\left(\frac{\lambda^{2}_{b}x_{b}}{|x_{b}|^{2}}+(b,0)\right).
\end{equation}
Multiplying the above identity by $|x|^{N-2+2\alpha}$ and letting
$|x|\rightarrow\infty,$ we find
$$\widetilde{u}(0)=\lim_{|x|\rightarrow\infty}|x|^{N-2+2\alpha}u(x)
=\lambda^{N-2+2\alpha}_{b}\lim_{|x|\rightarrow\infty}
\left(\frac{|x|}{|x_{b}|}\right)^{N-2+2\alpha}
u\left(\frac{\lambda^{2}_{b}x_{b}}{|x_{b}|^{2}}+(b,0)\right)=\lambda^{N-2+2\alpha}_{b}
u(b,0),$$ and using $u(0)=\widetilde{u}(0)$, we get
\begin{equation}\label{vxcewe}
\lambda^{N-2+2\alpha}_{b}=\frac{u(0)}{u(b,0)}.
\end{equation}
From $u=\widetilde{u}$ and (\ref{xcretyu}), we have
$\frac{1}{|x|^{N-2+2\alpha}}u\left(\frac{x}{|x|^{2}}\right)
=\left(\frac{\lambda_{b}}{|x_{b}|}\right)^{N-2+2\alpha}
u\left(\frac{\lambda^{2}_{b}x_{b}}{|x_{b}|^{2}}+(b,0)\right).$ Let
$f(x')=u(x',0),$ by Proposition \ref{tqaaa3w}, we know that $f\in
C^{1,\gamma}(\mathbb{R}^{N-1}).$ Now setting $x_{N}=0$ in the last
identity and using (\ref{vxcewe}), we obtain
$$\frac{1}{|x'|^{N-2+2\alpha}}f\left(\frac{x'}{|x'|^{2}}\right)
=\left(\frac{u(0)}{u(b,0)}\right)\frac{1}{|x'-b|^{N-2+2\alpha}}
f\left(\frac{\lambda^{2}_{b}(x'-b)}{|x'-b|^{2}}+b\right).$$ Then
as the proof of Corollary 2.8 of \cite{MM}, we can get that
$f(b)=f(x'_{0})(1+|b-x'_{0}|^{2})^{-(N-2+2\alpha)/2}$ for some
fixed $x'_{0}\in\mathbb{R}^{N-1}.$ By the arbitrariness of $b,$ we
have $u(x',0)=u(x'_{0},0)(1+|x'-x'_{0}|^{2})^{-(N-2+2\alpha)/2},$
$\forall x'\in\mathbb{R}^{N-1}.$\hfill$\Box$

\begin{corollary}\label{88eubbcv00ok} Suppose that $\alpha>1/2$ and $0\leq s<2.$
Let $u\in X_{\alpha,loc}(\mathbb{R}^{N})$ be a positive solution
of equation (\ref{yyyhsf2}). Then there exist $\lambda>0$ and
$x'_{0}\in \mathbb{R}^{N-1}$ such that
$u(x',0)=u(x'_{0},0)(1+\lambda^{2}|
x'-x'_{0}|^{2})^{-(N-2+2\alpha)/2},$ $\forall x'\in
\mathbb{R}^{N-1}$.
\end{corollary}
\noindent{\bf Proof.} By Theorem \ref{yqwe12}, there exists
$\mu>0$ such that $u_{\mu}\equiv u.$ Let
$v(x)=\mu^{N-2+2\alpha}u(\mu^{2} x)$, $x\in \mathbb{R}^{N}$. Then
$v$ is a solution of equation (\ref{yyyhsf2}) satisfying
$v(x)=|x|^{-(N-2+2\alpha)}v(x/|x|^{2})$, $x\in
\mathbb{R}^{N}\setminus \{0\}.$ By Theorem \ref{bz2wdf}, we get
that $v(x',0)=v(a,0)(1+|x'-a|^{2})^{-(N-2+2\alpha)/2}$, $\forall
x'\in \mathbb{R}^{N-1}$ for some $a\in \mathbb{R}^{N-1}$. By
$v(x)=\mu^{N-2+2\alpha}u(\mu^{2} x)$, we get that
$u(x',0)=u(x'_{0},0)(1+\lambda^{2}|
x'-x'_{0}|^{2})^{-(N-2+2\alpha)/2},$ $x'\in \mathbb{R}^{N-1}$,
with $\lambda=1/\mu^{2}$ and $x'_{0}=\mu^{2} a.$ \hfill$\Box$

\medskip

\begin{Theorem}\label{trgt88uwhs55r} Suppose that $\alpha>1/2$ and $0\leq s<2.$
Let $U_{\alpha,s}$ be a  positive solution of equation
(\ref{yyyhsf2}). Then $u$ is a positive solution of equation
(\ref{yyyhsf2}) if and only if
$u(x',x_{N})=\lambda^{\frac{N-2+2\alpha}{2}}U_{\alpha,s}(\lambda
x'+x'_{0},\lambda x_{N})$ for some $\lambda>0$ and $x'_{0}\in
\mathbb{R}^{N-1}$.
\end{Theorem}
\noindent{\bf Proof.} By Corollary \ref{88eubbcv00ok}, there exist
$\eta>0$ and $a\in \mathbb{R}^{N-1}$ such that
$U_{\alpha,s}(x',0)=U_{\alpha,s}(a,0)(1+\eta^{2}|x'-a|^{2})^{-(N-2+2\alpha)/2}.$
 From the proof of Proposition \ref{yhadsd}, we
know that $x^{\alpha}_{N}U_{\alpha,s}(x)$ and
$x^{\alpha}_{N}u(x)$,  $x\in \mathbb{R}^{N}_{+}$ are solutions of
equation (\ref{ooo889}). From Proposition 5.13 of \cite{CL}, we
deduce that
$W_{\alpha,s}(x)=x^{\frac{N-2+2\alpha}{2}}_{N}U_{\alpha,s}(x)$ and
$v(x)=x^{\frac{N-2+2\alpha}{2}}_{N}u(x)$, $x\in
\mathbb{R}^{N}_{+}$ are solutions of equation (\ref{ttfgbbcv456})
with $K\equiv1.$ By Remark \ref{ccv5}, we know that $U_{\alpha,s}$
and $u$ lie in the space $X_{\alpha}(\mathbb{R}^{N})$. Then it is
easy to verify that $W_{\alpha,s}$ and $v$ lie in the space $
H^{1}(\mathbb{H}),$ where $\mathbb{H}=(\mathbb{R}^{N}_{+},
\frac{dx^{2}}{x^{2}_{N}})$ is the $N-$dimensional hyperbolic space
and $H^{1}(\mathbb{H})$ is the Hilbert space
$\overline{C^{\infty}_{c}(\mathbb{H})}^{||\cdot||}$ with norm
$||u||=(\int_{\mathbb{\mathbb{H}}}|\nabla_{\mathbb{H}}u|^{2}dV_{\mathbb{H}})^{\frac{1}{2}}
=(\int_{\mathbb{R}^{N}_{+}}x^{2-N}_{N}|\nabla
u|^{2}dx)^{\frac{1}{2}}.$ From \cite{MS},  up to an isometric
transform of $\mathbb{H}$, the positive solution of equation
(\ref{ttfgbbcv456}) which lies in $H^{1}(\mathbb{H})$ is unique.
And  from page 116 of \cite{R}, we know that the isometric
transforms of $\mathbb{H}$ are those M\"obius transforms of
$\mathbb{R}^{N}$ which leave $\mathbb{R}^{N}_{+}$ invariant. Thus
there exist $\lambda>0$ and $x'_{0}\in \mathbb{R}^{N-1}$ such that
$v(x)=W_{\alpha,s}(\lambda x'+x'_{0},\lambda x_{N})$ for any
$x=(x',x_{N})\in \mathbb{R}^{N}_{+}$.  It follows that
$u(x',0)=\lambda^{\frac{N-2+2\alpha}{2}}U_{\alpha,s}(\lambda
x'+x'_{0},0),$ $\forall x'\in \mathbb{R}^{N-1}.$ Let
$\overline{W}_{\alpha,s}(x)=(-x_{N})^{\frac{N-2+2\alpha}{2}}U_{\alpha,s}(x)$
and $\overline{v}(x)=(-x_{N})^{\frac{N-2+2\alpha}{2}}u(x)$, $x\in
\mathbb{R}^{N}_{-}$. Using the same argument, we deduce that there
exist $\mu>0$ and $\overline{x}'_{0}\in \mathbb{R}^{N-1}$ such
that $\overline{v}(x)=\overline{W}_{\alpha,s}(\mu
x'+\overline{x}'_{0},\mu x_{N})$ for any $x=(x',x_{N})\in
\mathbb{R}^{N}_{-}$. In particular, we have
$u(x',0)=\mu^{\frac{N-2+2\alpha}{2}}U_{\alpha,s}(\mu
x'+\overline{x}'_{0},0). $ Thus
$\lambda^{\frac{N-2+2\alpha}{2}}U_{\alpha,s}(\lambda
x'+x'_{0},0)=\mu^{\frac{N-2+2\alpha}{2}}U_{\alpha,s}(\mu
x'+\overline{x}'_{0},0)$ for any $ x'\in \mathbb{R}^{N-1}.$ We
obtain $x'_{0}=\overline{x}'_{0}$ and $\lambda=\mu.$ Thus
$u(x',x_{N})=\lambda^{\frac{N-2+2\alpha}{2}}U_{\alpha,s}(\lambda
x'+x'_{0},\lambda x_{N})$ for some $\lambda>0$ and $x'_{0}\in
\mathbb{R}^{N-1}$. \hfill$\Box$

 \begin{Theorem}\label{ggdrefddf}
 If $u\in X_{1,loc}(\mathbb{R}^{N})$ is a positive solution of equation (\ref{yyyhsf2})
with $\alpha=1$ and $s=1+2/N,$ i.e., $u$ is a positive solution of
equation
 \begin{eqnarray}\label{tegdcc}
-div(|x_{N}|^{2}\nabla u)=|x_{N}|u^{\frac{N+2}{N}},
\end{eqnarray}
then there exist $\lambda>0$ and $\zeta\in \mathbb{R}^{N-1}$ such
that $u(x',x_{N})=\lambda^{\frac{N}{2}}U(\lambda x'+\zeta,\lambda
x_{N})$, where
$U(x',x_{N})=\left(\frac{2N}{(1+|x_{N}|)^{2}+|x'|^{2}}\right)^{\frac{N}{2}}.$
  Furthermore, taking derivatives with
respect to the parameters $\lambda$ and $\zeta$ at $\lambda=1$ and
$\zeta=0,$ we get $N$ functions $V_{1},\cdots,V_{N}.$ These
functions are solutions to the linearized equation
\begin{eqnarray}\label{jdgcvvv5}
-div(|x_{N}|^{2}\nabla v)=\frac{N+2}{N}|x_{N}|U^{\frac{2}{N}}v\
\mbox{in}\ \mathbb{R}^{N},\quad v\in X_{1}(\mathbb{R}^{N}),
\end{eqnarray}
and any solution of (\ref{jdgcvvv5}) can be the linear combination
of $V_{1},\cdots,V_{N}.$
 \end{Theorem}
\noindent{\bf Proof.} If $u$ is a positive  solution of equation
(\ref{tegdcc}), then by Remark \ref{ccv5}, we know that $u\in
X_{1}(\mathbb{R}^{N})$. Then $v=x_{N}u\in
\mathcal{D}^{1,2}_{0}(\mathbb{R}^{N}_{+})$ (see (\ref{oooo91s}))
and it is a positive solution of equation $ -\triangle
v=\frac{v^{(N+2)/N}}{x^{1+2/N}_{N}} \ \mbox{in}\
\mathbb{R}^{N}_{+}$ (see ( \ref{ooo889})).  From Proposition 5.13
of \cite{CL}, we know that $v_{1}(x)=x_{N}^{\frac{N-2}{2}}
v(x)=x^{\frac{N}{2}}_{N}u$ is a solution of the  equation
$-\triangle_{\mathbb{H}}v_{1}=\frac{N(N-2)}{4}v_{1}+v^{\frac{N+2}{N}}_{1}$
which satisfies $v_{1}\in H^{1}(\mathbb{H}),$ where
$\mathbb{H}=(\mathbb{R}^{N}_{+}, \frac{dx^{2}}{x^{2}_{N}})$ is the
$N-$dimensional hyperbolic space and $H^{1}(\mathbb{H})$ is the
Sobolev space defined in the proof of Theorem \ref{trgt88uwhs55r}.
Let $\mathbb{R}^{N+2}=\mathbb{R}^{N-1}\times \mathbb{R}^{3}$ and
$z=(x,y)$, $x\in \mathbb{R}^{N-1},$ $y\in\mathbb{R}^{3}.$ By
\cite[Lemma 2.1]{cfms}, we know that
$u_{1}(x,y)=|y|^{-\frac{N}{2}}v_{1}(x,|y|)=u(x,|y|)$ is a solution
of equation $-\triangle u_{1}=\frac{u^{\frac{N+2}{N}}_{1}}{|y|}$
with $u_{1}\in \mathcal{D}^{1,2}_{0}(\mathbb{R}^{N+2}).$ By
\cite[Theorem 1.1]{MFS},  Up to dilations and translations in $x$,
this equation has unique solution
$U_{1}(x,y)=\left(\frac{2N}{(1+|y|)^{2}+|x|^{2}}\right)^{\frac{N}{2}}.$
Therefore,   up to dilations and translations in $x'$, equation
(\ref{tegdcc}) has a unique positive solution
$U(x',x_{N})=\left(\frac{2N}{(1+|x_{N}|)^{2}+|x'|^{2}}\right)^{\frac{N}{2}}.$
By \cite[Theorem 3.1]{cfms}, taking derivatives with respect to
the parameters $\lambda$ and $\zeta$ at $\lambda=1$ and $\zeta=0$
to $\lambda^{(N-2)/2}U_{1}(\lambda x+\zeta,\lambda y)$, we get $N$
functions. These functions are solutions to the linearized
equation at $U_{1}$
\begin{eqnarray}\label{jff222v5}
-\triangle v=\frac{N+2}{N}\frac{U^{\frac{2}{N}}_{1}}{|y|}v \quad
\mbox{in}\ \mathbb{R}^{N+2},\quad v\in
\mathcal{D}^{1,2}_{0}(\mathbb{R}^{N+2})
\end{eqnarray}
and any solution of (\ref{jff222v5}) can be the linear combination
of the $N$ functions. Thus taking derivatives with respect to the
parameters $\lambda$ and $\zeta$ at $\lambda=1$ and $\zeta=0$ to
$\lambda^{N/2}U(\lambda x'+\zeta,\lambda x_{N})$,  we get $N$
functions $V_{1},\cdots,V_{N}.$ These functions are solutions to
the linearized equation (\ref{jdgcvvv5})   and  any solution of
(\ref{jdgcvvv5}) can be the linear combination of
$V_{1},\cdots,V_{N}.$ \hfill$\Box$

\section{Some variational  identities}\label{88ujzzx}
\hspace*{\parindent} In this section, we derive some variational
identities for solutions of equation (\ref{yyyhsf}). As a
consequence, some nonexistence results for solutions of  equation
(\ref{yyyhsf}) are obtained.
\begin{Theorem}\label{bvcxq}
If $K\in C^{1}(\overline{B_{\varsigma}(0))}$   and $u\in
X_{\alpha,loc}(B_{\varsigma}(0))$ is a weak   solution of equation
(\ref{yyyhsf}) in $B_{\varsigma}(0),$ then for any
$0<\sigma<\varsigma,$  the following identity holds
\begin{eqnarray}\label{xzc4423}
&&\frac{1}{2^{*}(s)}\int_{B_{\sigma}(0)}(x\cdot\nabla
K)\cdot|x_{N}|^{\alpha\cdot2^{*}(s)-s}|u|^{2^{*}(s)}-
\frac{1}{2^{*}(s)}\int_{\partial
B_{\sigma}(0)}(x\cdot\mathbf{n})\cdot K(x)
|x_{N}|^{\alpha\cdot2^{*}(s)-s}|u|^{2^{*}(s)}\nonumber\\
&=& \int_{\partial B_{\sigma}(0)}B(\sigma,x,u,\nabla u),
\end{eqnarray}
where $B(\sigma,x,u,\nabla
u)=\frac{N-2+2\alpha}{2}|x_{N}|^{2\alpha}\cdot u\frac{\partial
u}{\partial\mathbf{n}}-\frac{\sigma}{2}|x_{N}|^{2\alpha}|\nabla
u|^{2}+\sigma|x_{N}|^{2\alpha}\left(\frac{\partial
u}{\partial\mathbf{ n}}\right)^{2}$ and
$\mathbf{n}=(n_{1},\cdots,n_{N})$ is the outer normal vector of
$\partial B_{\sigma}(0)$, i.e., $\mathbf{n}=x/|x|,$
$n_{i}=x_{i}/|x|,$ $1\leq i\leq N.$
\end{Theorem}
\noindent{\bf Proof.} For $0<\epsilon<\sigma,$ let
$\Omega^{+}_{\epsilon,\sigma}=B_{\sigma}(0)\cap\{x_{N}>\epsilon\}.$
Multiplying left hand side of equation (\ref{yyyhsf}) by
$x\cdot\nabla u$ and integrating in
$\Omega^{+}_{\epsilon,\sigma},$ we obtain by divergence theorem
that
\begin{eqnarray}\label{vxc423w}
&&-\int_{\Omega^{+}_{\epsilon,\sigma}}div(|x_{N}|^{2\alpha}\nabla
u)(x\cdot\nabla u) \nonumber\\
&=&-\int_{\partial\Omega^{+}_{\epsilon,\sigma}}
|x_{N}|^{2\alpha}(\nabla u\cdot\mathbf{n})(x\cdot\nabla u)+
\int_{\Omega^{+}_{\epsilon,\sigma}}|x_{N}|^{2\alpha}\nabla
u\cdot\nabla (x\cdot\nabla u),
\end{eqnarray}
where $\mathbf{n}=(n_{1},\cdots,n_{N})$ is the outer normal vector
of $\Omega^{+}_{\epsilon,\sigma}.$ We have
\begin{eqnarray}\label{cvdfe1}
\int_{\Omega^{+}_{\epsilon,\sigma}}|x_{N}|^{2\alpha}\nabla
u\cdot\nabla (x\cdot\nabla
u)=\int_{\Omega^{+}_{\epsilon,\sigma}}|x_{N}|^{2\alpha}|\nabla
u|^{2} +\sum^{N}_{i=1}\sum^{N}_{j=1}
\int_{\Omega^{+}_{\epsilon,\sigma}}|x_{N}|^{2\alpha}x_{j}\frac{\partial
u}{\partial x_{i}}\frac{\partial^{2} u}{\partial x_{i}\partial
x_{j}}.
\end{eqnarray}
Through integrating by part, we get that
\begin{eqnarray}\label{b6667}
&&\sum^{N}_{i=1}\sum^{N}_{j=1}
\int_{\Omega^{+}_{\epsilon,\sigma}}|x_{N}|^{2\alpha}x_{j}\frac{\partial
u}{\partial x_{i}}\frac{\partial^{2} u}{\partial x_{i}\partial
x_{j}}\nonumber\\
&=&\sum^{N}_{i=1}\sum^{N}_{j=1}
\int_{\partial\Omega^{+}_{\epsilon,\sigma}}|x_{N}|^{2\alpha}(n_{j}\cdot
x_{j})\cdot\left(\frac{\partial u}{\partial x_{i}}\right)^{2}-
\sum^{N}_{i=1}\sum^{N}_{j=1} \int_{\Omega^{+}_{\epsilon,\sigma}}
\frac{\partial u}{\partial x_{i}}\cdot\frac{\partial}{\partial
x_{j}}\left(|x_{N}|^{2\alpha} x_{j}\frac{\partial u}{\partial
x_{i}}\right)\nonumber\\
&=&\int_{\partial\Omega^{+}_{\epsilon,\sigma}}|x_{N}|^{2\alpha}(\mathbf{n}\cdot
x)\cdot|\nabla u|^{2}\nonumber\\
&&-N\int_{\Omega^{+}_{\epsilon,\sigma}}|x_{N}|^{2\alpha}|\nabla
u|^{2}-2\alpha\int_{\Omega^{+}_{\epsilon,\sigma}}|x_{N}|^{2\alpha}
|\nabla u|^{2}-\sum^{N}_{i=1}\sum^{N}_{j=1}
\int_{\Omega^{+}_{\epsilon,\sigma}}|x_{N}|^{2\alpha}x_{j}\frac{\partial
u}{\partial x_{i}}\frac{\partial^{2} u}{\partial x_{i}\partial
x_{j}}.\nonumber
\end{eqnarray}
It follows that
\begin{equation}\label{vx4e}
\sum^{N}_{i=1}\sum^{N}_{j=1}
\int_{\Omega^{+}_{\epsilon,\sigma}}|x_{N}|^{2\alpha}x_{j}\frac{\partial
u}{\partial x_{i}}\frac{\partial^{2} u}{\partial x_{i}\partial
x_{j}}=\frac{1}{2}\int_{\partial\Omega^{+}_{\epsilon,\sigma}}|x_{N}|^{2\alpha}(\mathbf{n}\cdot
x)\cdot|\nabla
u|^{2}-\frac{N+2\alpha}{2}\int_{\Omega^{+}_{\epsilon,\sigma}}|x_{N}|^{2\alpha}
|\nabla u|^{2}.
\end{equation}
By $(\ref{vxc423w})-(\ref{vx4e}),$ we obtain
\begin{eqnarray}\label{xc445}
-\int_{\Omega^{+}_{\epsilon,\sigma}}div(|x_{N}|^{2\alpha}\nabla
u)(x\cdot\nabla u) &=&-\int_{\partial\Omega^{+}_{\epsilon,\sigma}}
|x_{N}|^{2\alpha}(\nabla u\cdot\mathbf{n})(x\cdot\nabla u)
+\frac{1}{2}\int_{\partial\Omega^{+}_{\epsilon,\sigma}}|x_{N}|^{2\alpha}(\mathbf{n}\cdot
x)\cdot|\nabla
u|^{2}\nonumber\\
&&-\frac{N-2+2\alpha}{2}\int_{\Omega^{+}_{\epsilon,\sigma}}|x_{N}|^{2\alpha}
|\nabla u|^{2}.
\end{eqnarray}
Multiplying right hand side of equation (\ref{yyyhsf}) by
$x\cdot\nabla u$ and integrating in
$\Omega^{+}_{\epsilon,\sigma},$ we obtain
\begin{eqnarray}\label{cvf8}
&&\int_{\Omega^{+}_{\epsilon,\sigma}}K(x)|x_{N}|^{\alpha\cdot
2^{*}(s)-s}|u|^{2^{*}(s)-2}u\cdot(x\cdot\nabla u)\nonumber\\
&=&\frac{1}{2^{*}(s)}\sum^{N}_{i=1}\int_{\Omega^{+}_{\epsilon,\sigma}}K(x)|x_{N}|^{\alpha\cdot
2^{*}(s)-s}x_{i}\frac{\partial}{\partial
x_{i}}\left(|u|^{2^{*}(s)}\right)\nonumber\\
&=&\frac{1}{2^{*}(s)}\left\{\int_{\partial\Omega^{+}_{\epsilon,\sigma}}K(x)|x_{N}|^{\alpha\cdot
2^{*}(s)-s}(x\cdot \mathbf{n})\cdot
|u|^{2^{*}(s)}-\sum^{N}_{i=1}\int_{\Omega^{+}_{\epsilon,\sigma}}
|u|^{2^{*}(s)}\frac{\partial}{\partial
x_{i}}\left(K(x)|x_{N}|^{\alpha\cdot
2^{*}(s)-s}x_{i}\right)\right\}\nonumber\\
&=&\frac{1}{2^{*}(s)}\int_{\partial\Omega^{+}_{\epsilon,\sigma}}K(x)|x_{N}|^{\alpha\cdot
2^{*}(s)-s}\cdot(x\cdot\mathbf{n})\cdot
|u|^{2^{*}(s)}-\frac{N-2+2\alpha}{2}\int_{\Omega^{+}_{\epsilon,\sigma}}K(x)|x_{N}|^{\alpha\cdot
2^{*}(s)-s}|u|^{2^{*}(s)}\nonumber\\
&&-\frac{1}{2^{*}(s)}\int_{\Omega^{+}_{\epsilon,\sigma}}(x\cdot\nabla
K(x))|x_{N}|^{\alpha\cdot 2^{*}(s)-s} |u|^{2^{*}(s)}.
\end{eqnarray}
By (\ref{xc445}), (\ref{cvf8}) and the fact that
$-\int_{\Omega^{+}_{\epsilon,\sigma}}div(|x_{N}|^{2\alpha}\nabla
u)(x\cdot\nabla
u)=\int_{\Omega^{+}_{\epsilon,\sigma}}K(x)|x_{N}|^{\alpha\cdot
2^{*}(s)-s}|u|^{2^{*}(s)-2}u\cdot(x\cdot\nabla u),$ we have
\begin{eqnarray}\label{zxcrrf}
&&-\int_{\partial\Omega^{+}_{\epsilon,\sigma}}
|x_{N}|^{2\alpha}(\nabla u\cdot\mathbf{n})(x\cdot\nabla u)
+\frac{1}{2}\int_{\partial\Omega^{+}_{\epsilon,\sigma}}|x_{N}|^{2\alpha}(\mathbf{n}\cdot
x)\cdot|\nabla
u|^{2}-\frac{N-2+2\alpha}{2}\int_{\Omega^{+}_{\epsilon,\sigma}}|x_{N}|^{2\alpha}
|\nabla u|^{2}\nonumber\\
&=&\frac{1}{2^{*}(s)}\int_{\partial\Omega^{+}_{\epsilon,\sigma}}K(x)|x_{N}|^{\alpha\cdot
2^{*}(s)-s}\cdot(x\cdot\mathbf{n})\cdot |u|^{2^{*}(s)}\nonumber\\
&&-\frac{N-2+2\alpha}{2}\int_{\Omega^{+}_{\epsilon,\sigma}}K(x)|x_{N}|^{\alpha\cdot
2^{*}(s)-s}|u|^{2^{*}(s)}
-\frac{1}{2^{*}(s)}\int_{\Omega^{+}_{\epsilon,\sigma}}(x\cdot\nabla
K(x))|x_{N}|^{\alpha\cdot 2^{*}(s)-s} |u|^{2^{*}(s)}.
\end{eqnarray}
Since
$-\int_{\Omega^{+}_{\epsilon,\sigma}}div(|x_{N}|^{2\alpha}\nabla
u)u=\int_{\Omega^{+}_{\epsilon,\sigma}}K(x)|x_{N}|^{\alpha\cdot
2^{*}(s)-s}|u|^{2^{*}(s)}$ and
$$-\int_{\Omega^{+}_{\epsilon,\sigma}}div(|x_{N}|^{2\alpha}\nabla
u)u=-\int_{\partial\Omega^{+}_{\epsilon,\sigma}}|x_{N}|^{2\alpha}(\mathbf{n}\cdot\nabla
u)\cdot
u+\int_{\Omega^{+}_{\epsilon,\sigma}}|x_{N}|^{2\alpha}|\nabla
u|^{2},$$ we have
\begin{eqnarray}\label{vxcer1}
-\int_{\partial\Omega^{+}_{\epsilon,\sigma}}|x_{N}|^{2\alpha}(\mathbf{n}\cdot\nabla
u)\cdot
u+\int_{\Omega^{+}_{\epsilon,\sigma}}|x_{N}|^{2\alpha}|\nabla
u|^{2}=\int_{\Omega^{+}_{\epsilon,\sigma}}K(x)|x_{N}|^{\alpha\cdot
2^{*}(s)-s}|u|^{2^{*}(s)}.
\end{eqnarray}
By (\ref{zxcrrf}) and (\ref{vxcer1}), we obtain
\begin{eqnarray}\label{czce45e}
&&\frac{1}{2^{*}(s)}\int_{\Omega^{+}_{\epsilon,\sigma}}(x\cdot\nabla
K(x))|x_{N}|^{\alpha\cdot 2^{*}(s)-s} |u|^{2^{*}(s)}
-\frac{1}{2^{*}(s)}\int_{\partial\Omega^{+}_{\epsilon,\sigma}}K(x)|x_{N}|^{\alpha\cdot
2^{*}(s)-s}\cdot(x\cdot\mathbf{n})\cdot |u|^{2^{*}(s)}\nonumber\\
&=&\frac{N-2+2\alpha}{2}
\int_{\partial\Omega^{+}_{\epsilon,\sigma}}|x_{N}|^{2\alpha}(\mathbf{n}\cdot\nabla
u)\cdot
u-\frac{1}{2}\int_{\partial\Omega^{+}_{\epsilon,\sigma}}|x_{N}|^{2\alpha}(\mathbf{n}\cdot
x)\cdot|\nabla
u|^{2}\nonumber\\
&&+\int_{\partial\Omega^{+}_{\epsilon,\sigma}}
|x_{N}|^{2\alpha}(\nabla u\cdot\mathbf{n})(x\cdot\nabla u).
\end{eqnarray}
Let
$\Omega^{-}_{\epsilon,\sigma}=B_{\sigma}(0)\cap\{x_{N}<-\epsilon\},$
$0<\epsilon <\sigma.$ We can get a similar identity like
(\ref{czce45e}) with $\Omega^{+}_{\epsilon,\sigma}$  replaced by
$\Omega^{-}_{\epsilon,\sigma}$. Adding these two identities, we
obtain
\begin{eqnarray}\label{czce45xe}
&&\frac{1}{2^{*}(s)}\int_{\Omega^{+}_{\epsilon,\sigma}\cup\Omega^{-}_{\epsilon,\sigma}}(x\cdot\nabla
K(x))|x_{N}|^{\alpha\cdot 2^{*}(s)-s}
|u|^{2^{*}(s)}\nonumber\\
&&-\frac{1}{2^{*}(s)}\int_{\partial\Omega^{+}_{\epsilon,\sigma}
\cup\partial\Omega^{-}_{\epsilon,\sigma}}K(x)|x_{N}|^{\alpha\cdot
2^{*}(s)-s}\cdot(x\cdot\mathbf{n})\cdot |u|^{2^{*}(s)}\nonumber\\
&=&\frac{N-2+2\alpha}{2}
\int_{\partial\Omega^{+}_{\epsilon,\sigma}\cup
\partial\Omega^{-}_{\epsilon,\sigma}}|x_{N}|^{2\alpha}(\mathbf{n}\cdot\nabla
u)\cdot u-\frac{1}{2}\int_{\partial\Omega^{+}_{\epsilon,\sigma}
\cup\partial\Omega^{-}_{\epsilon,\sigma}}|x_{N}|^{2\alpha}(\mathbf{n}\cdot
x)\cdot|\nabla
u|^{2}\nonumber\\
&&+\int_{\partial\Omega^{+}_{\epsilon,\sigma}
\cup\partial\Omega^{-}_{\epsilon,\sigma}} |x_{N}|^{2\alpha}(\nabla
u\cdot\mathbf{n})(x\cdot\nabla u),
\end{eqnarray}
where $\mathbf{n}=(n_{1},\cdots,n_{N})$ is the outer normal vector
of $\Omega^{+}_{\epsilon,\sigma}\cup\Omega^{-}_{\epsilon,\sigma}$.
 Since
$\partial\Omega^{+}_{\epsilon,\sigma}\cup\partial\Omega^{-}_{\epsilon,\sigma}
=(\{|x_{N}|>\epsilon\}\cap\partial B_{\sigma}(0))\cup
(\{|x_{N}|=\epsilon\}\cap B_{\sigma}(0))$, we get that  $
\partial\Omega^{+}_{\epsilon,\sigma}\cup\partial\Omega^{-}_{\epsilon,\sigma}
\rightarrow\partial B_{\sigma}(0)\cup ((\{x_{N}=0\}\cap
B_{\sigma}(0)))$ as $\epsilon\rightarrow 0.$ Moreover,
$\Omega^{+}_{\epsilon,\sigma}\cup\Omega^{-}_{\epsilon,\sigma}\rightarrow
B_{\sigma}(0)$ as $\epsilon\rightarrow 0.$ By  Proposition
\ref{fff1fc}, we get that as $\epsilon\rightarrow 0,$
\begin{eqnarray}\label{hhde}
\int_{\{x_{N}=\pm\epsilon\}\cap
B_{\sigma}(0)}|x_{N}|^{2\alpha}(\mathbf{n}\cdot\nabla u)\cdot u
=\mp\int_{\{x_{N}=\pm\epsilon\}\cap
B_{\sigma}(0)}\epsilon^{2\alpha}\frac{\partial u}{\partial
x_{N}}(x',\epsilon)\cdot u\rightarrow 0.
\end{eqnarray}
Furthermore, by Proposition \ref{tqaaa3w} and Proposition
\ref{fff1fc}, we deduce that as $\epsilon\rightarrow 0,$
\begin{eqnarray}\label{ccffxcd}
&&\left|\int_{\{x_{N}=\pm\epsilon\}\cap
B_{\sigma}(0)}|x_{N}|^{2\alpha}(\mathbf{n}\cdot x)\cdot|\nabla
u|^{2}\right|\nonumber\\
&\leq&C\int_{\{x_{N}=\pm\epsilon\}\cap
B_{\sigma}(0)}\epsilon^{2\alpha+1}\left(\left|\frac{\partial
u}{\partial x_{N}}(x',\epsilon)\right|^{2}+\left|\frac{\partial
u}{\partial x'}(x',\epsilon)\right|^{2}\right)\rightarrow 0
\end{eqnarray}
In a similar manner, we have
\begin{eqnarray}\label{ccff}
\lim_{\epsilon\rightarrow 0}\int_{\{x_{N}=\pm\epsilon\}\cap
B_{\sigma}(0)} |x_{N}|^{2\alpha}(\nabla
u\cdot\mathbf{n})(x\cdot\nabla u)=0.
\end{eqnarray}
Letting $\epsilon\rightarrow 0$ in (\ref{czce45xe}), by
$(\ref{hhde})-(\ref{ccff}),$ we obtain the desired result of this
theorem.\hfill$\Box$

\begin{corollary}\label{9981}
Suppose  $K, \widehat{K}\in C^{1}(\mathbb{R}^{N})\cap
L^{\infty}(\mathbb{R}^{N})$, where $\widehat{K}(x)=K(x/|x|).$ If
$u\in X_{\alpha}(\mathbb{R}^{N})$ is a weak  solution of equation
(\ref{yyyhsf}) in $\mathbb{R}^{N}$, then  $
\int_{\mathbb{R}^{N}}(x\cdot\nabla
K)\cdot|x_{N}|^{\alpha\cdot2^{*}(s)-s}|u|^{2^{*}(s)}(x)=0. $
\end{corollary}
\noindent{\bf Proof.} Let
$\widetilde{u}(y)=\frac{1}{|y|^{N-2+2\alpha}}u\left(\frac{y}{|y|^{2}}\right).$
Then
\begin{eqnarray}\label{gdbcg66dydgg}
&&\int_{\mathbb{R}^N}|y_N|^{2\alpha}|\nabla_y
\widetilde{u}|^2\nonumber\\
&=&\int_{\mathbb{R}^N}|y_N|^{2\alpha}\left|\nabla_y
\left(\frac{1}{|y|^{N-2+2\alpha}}u\left(\frac{y}{|y|^{2}}\right)\right)\right|^2\nonumber\\
&=&\int_{\mathbb{R}^N}|y_N|^{2\alpha}\left|
\frac{N-2+2\alpha}{|y|^{N-1+2\alpha}}\frac{y}{|y|}u\left(\frac{y}{|y|^{2}}\right)+
\frac{1}{|y|^{N-2+2\alpha}}\nabla_y\left(u\left(\frac{y}{|y|^{2}}\right)\right)\right|^2\nonumber\\
&\leq& 2\int_{\mathbb{R}^N}|y_N|^{2\alpha}\left|
\frac{1}{|y|^{N-1+2\alpha}}\frac{y}{|y|}u\left(\frac{y}{|y|^{2}}\right)\right|^2+
2\int_{\mathbb{R}^N}|y_N|^{2\alpha}
\left|\frac{1}{|y|^{N-2+2\alpha}}\nabla_y\left(u\left(\frac{y}{|y|^{2}}\right)\right)\right|^2.\nonumber\\
\end{eqnarray}
Using the transform $y=x/|x|^{2}$ ( the Jacobian of this transform
is $|x|^{-2N}$ ), we get that
\begin{eqnarray}\label{bfb99fijfjj}
\int_{\mathbb{R}^N}|y_N|^{2\alpha}\left|
\frac{1}{|y|^{N-1+2\alpha}}\frac{y}{|y|}u\left(\frac{y}{|y|^{2}}\right)\right|^2
&=&\int_{\mathbb{R}^N}|y_N|^{2\alpha}
\frac{1}{|y|^{2N-2+4\alpha}}u^2\left(\frac{y}{|y|^{2}}\right)\nonumber\\
&=&\int_{\mathbb{R}^N}|x_N|^{2\alpha}|x|^{-2}u^2(x)dx.
\end{eqnarray}
From the proof of Theorem \ref{yyyrtrtww}, we get that
$v(x):=|x_N|^\alpha u(x),$ $x\in \mathbb{R}^N$ satisfies that
$v\in H^1(\mathbb{R}^N)$. Then by the Hardy inequality, we deduce
that $\int_{\mathbb{R}^N}\frac{v^2}{|x|^2}<\infty.$ Thus,
\begin{eqnarray}\label{hhfhfgg88ey}
\int_{\mathbb{R}^N}|x_N|^{2\alpha}|x|^{-2}u^2(x)dx=\int_{\mathbb{R}^N}\frac{v^2}{|x|^2}<\infty.
\end{eqnarray}
By (\ref{bfb99fijfjj}) and (\ref{hhfhfgg88ey}), we get that
\begin{eqnarray}\label{hhfjjjjr7746ttg}
\int_{\mathbb{R}^N}|y_N|^{2\alpha}\left|
\frac{1}{|y|^{N-1+2\alpha}}\frac{y}{|y|}u\left(\frac{y}{|y|^{2}}\right)\right|^2<\infty.
\end{eqnarray}

 Note that
\begin{eqnarray}\label{gfb9fijjfgsvz}
\nabla_y\left(u\left(\frac{y}{|y|^{2}}\right)\right)=(\nabla_x
u)\left(\frac{y}{|y|^2}\right)\cdot A(y),
\end{eqnarray}
where
$A(y)=\left(\frac{\delta_{i,j}}{|y|^2}-\frac{2y_iy_j}{|y|^4}\right)_{N\times
N}$ which satisfies that
\begin{eqnarray}\label{hvbnb88gugyg}
A(y)\cdot\frac{y}{|y|}=-\frac{y}{|y|^3},\
A(y)A(y)^T=\frac{1}{|y|^4}I.
\end{eqnarray}
It follows that
\begin{eqnarray}\label{gfbvhcydttegd}
\left|\nabla_y\left(u\left(\frac{y}{|y|^{2}}\right)\right)\right|^2&=&(\nabla_x
u)\left(\frac{y}{|y|^2}\right)\cdot A(y)A(y)^T(\nabla_x
u)\left(\frac{y}{|y|^2}\right)^T\nonumber\\
&=&\frac{1}{|y|^4}\left|(\nabla_x
u)\left(\frac{y}{|y|^{2}}\right)\right|^2.
\end{eqnarray}
By (\ref{gfbvhcydttegd}) and using the transform $y=x/|x|^{2}$, we
get that
\begin{eqnarray}\label{jjgnbjghughg09i}
\int_{\mathbb{R}^N}|y_N|^{2\alpha}
\left|\frac{1}{|y|^{N-2+2\alpha}}\nabla_y\left(u\left(\frac{y}{|y|^{2}}\right)\right)\right|^2
&=&\int_{\mathbb{R}^N}\frac{|y_N|^{2\alpha}}{|y|^{2N+4\alpha}}\left|(\nabla_x
u)\left(\frac{y}{|y|^{2}}\right)\right|^2\nonumber\\
&=&\int_{\mathbb{R}^N}|x_N|^{2\alpha}|\nabla_x u|^2<\infty.
\end{eqnarray}
By (\ref{gdbcg66dydgg}), (\ref{hhfjjjjr7746ttg}) and
(\ref{jjgnbjghughg09i}), we get that
$\int_{\mathbb{R}^N}|y_N|^{2\alpha}|\nabla_y
\widetilde{u}|^2<\infty.$ It follows that $\widetilde{u}\in
X_\alpha.$ Thus by Theorem \ref{tererd}, we get that
$\widetilde{u}\in L^\infty(\mathbb{R}^N).$

Since
\begin{eqnarray}\label{jianshezhe8dydy}
\nabla_y \widetilde{u}(y)=-
\frac{N-2+2\alpha}{|y|^{N-1+2\alpha}}\frac{y}{|y|}u\left(\frac{y}{|y|^{2}}\right)+
\frac{1}{|y|^{N-2+2\alpha}}\nabla_y\left(u\left(\frac{y}{|y|^{2}}\right)\right),
\end{eqnarray}
by (\ref{gfbvhcydttegd}), we get that
\begin{eqnarray}\label{hnavcbbf9fifj}
&&\sigma\int_{\partial B_\sigma(0)}|y_{N}|^{2\alpha}|\nabla_y
\widetilde{u}|^{2}\nonumber\\
&\geq&\frac{\sigma}{2}\int_{\partial
B_\sigma(0)}\frac{|y_N|^{2\alpha}}{|y|^{2N+4\alpha}}\left|(\nabla_x
u)\left(\frac{y}{|y|^{2}}\right)\right|^2
-\sigma(N-2+2\alpha)^2\int_{\partial B_\sigma(0)}
\frac{|y_N|^{2\alpha}}{|y|^{2N-2+4\alpha}}u^2\left(\frac{y}{|y|^{2}}\right)\nonumber\\
&=&\frac{\sigma}{2\sigma^{2N+4\alpha}}\int_{\partial
B_\sigma(0)}|y_N|^{2\alpha}\left|(\nabla_x
u)\left(\frac{y}{|y|^{2}}\right)\right|^2-\sigma(N-2+2\alpha)^2\int_{\partial
B_\sigma(0)} \frac{|y_N|^{2\alpha}}{|y|^{2}}\widetilde{u}^2(y).
\end{eqnarray}
Using the transform $y=x/|x|^{2}$, we get that
\begin{eqnarray}\label{hfbvh77fyfh}
\frac{\sigma}{\sigma^{2N+4\alpha}}\int_{\partial
B_\sigma(0)}|y_N|^{2\alpha}\left|(\nabla_x
u)\left(\frac{y}{|y|^{2}}\right)\right|^2
&=&\frac{\sigma}{\sigma^{2N+4\alpha}}\int_{\partial
B_{1/\sigma}(0)}\frac{|x_N|^{2\alpha}}{|x|^{4\alpha}}|(\nabla_x
u)(x)|^2\cdot\sigma^{2(N-1)}\nonumber\\
&=&\sigma^{-1}\int_{\partial
B_{1/\sigma}(0)}|x_N|^{2\alpha}|(\nabla_x u)(x)|^2.
\end{eqnarray}
Since $\widetilde{u}\in X_\alpha$, as (\ref{hhfhfgg88ey}), we can
get that
\begin{eqnarray}\label{hfhb99fijfhuuw}
\lim_{\sigma\rightarrow 0}\sigma\int_{\partial
B_\sigma(0)}|y_{N}|^{2\alpha}|\nabla_y \widetilde{u}|^{2}=0,\
\lim_{\sigma\rightarrow 0}\sigma(N-2+2\alpha)^2\int_{\partial
B_\sigma(0)}\frac{ |y_N|^{2\alpha}}{|y|^2}\widetilde{u}^2(y)=0.
\end{eqnarray}
Then by $(\ref{hnavcbbf9fifj})-(\ref{hfhb99fijfhuuw})$, we get
that
\begin{eqnarray}\label{hfhbcgdteyye00}
\lim_{\sigma\rightarrow 0}\sigma^{-1}\int_{\partial
B_{1/\sigma}(0)}|x_N|^{2\alpha}|(\nabla_x u)(x)|^2=0.
\end{eqnarray}
It follows that
\begin{eqnarray}\label{hfb77fyfgggg}
\lim_{\sigma\rightarrow 0}\sigma^{-1}\int_{\partial
B_{1/\sigma}(0)}|x_{N}|^{2\alpha}\left(\frac{\partial
u}{\partial\mathbf{ n}}\right)^{2}=0.
\end{eqnarray}
By (\ref{jianshezhe8dydy}), (\ref{gfb9fijjfgsvz}) and
(\ref{hvbnb88gugyg}), we get
\begin{eqnarray}\label{hhfbv99didhhh}
\nabla_y \widetilde{u}\cdot\frac{y}{|y|}&=&-
\frac{N-2+2\alpha}{|y|^{N-1+2\alpha}}u\left(\frac{y}{|y|^{2}}\right)
+\frac{1}{|y|^{N-2+2\alpha}}(\nabla_x
u)\left(\frac{y}{|y|^2}\right)\cdot
A(y)\cdot\frac{y}{|y|}\nonumber\\
&=&-
\frac{N-2+2\alpha}{|y|^{N-1+2\alpha}}u\left(\frac{y}{|y|^{2}}\right)
-\frac{1}{|y|^{N+2\alpha}}(\nabla_x
u)\left(\frac{y}{|y|}\right)\cdot
\frac{y}{|y|}\nonumber\\
&=&-
(N-2+2\alpha)\frac{\widetilde{u}(y)}{|y|}-\frac{1}{|y|^{N+2\alpha}}(\nabla_x
u)\left(\frac{y}{|y|}\right)\cdot
\frac{y}{|y|}.\nonumber\end{eqnarray} It follows that
\begin{eqnarray}\label{ggfbjjfufhhhf77qe}
&&\int_{\partial B_\sigma(0)}|y_N|^{2\alpha}\widetilde{u}\nabla_y
\widetilde{u}\cdot\frac{y}{|y|}\nonumber\\
&=&-(N-2+2\alpha)\int_{\partial
B_\sigma(0)}\frac{|y_N|^{2\alpha}}{|y|} \widetilde{u}^2(y)
-\int_{\partial
B_\sigma(0)}\frac{|y_N|^{2\alpha}}{|y|^{2N-2+4\alpha}}u\left(\frac{y}{|y|}\right)\left((\nabla_x
u)\left(\frac{y}{|y|}\right)\cdot \frac{y}{|y|}\right).\nonumber\\
\end{eqnarray}
Using the transform $y=x/|x|^{2}$, we get that
\begin{eqnarray}\label{iidhdhgtdrddfd}
&&\int_{\partial
B_\sigma(0)}\frac{|y_N|^{2\alpha}}{|y|^{2N-4+4\alpha}}u\left(\frac{y}{|y|}\right)\left((\nabla_x
u)\left(\frac{y}{|y|}\right)\cdot \frac{y}{|y|}\right)\nonumber\\
&=&\frac{1}{\sigma^{2N-2+4\alpha}}\int_{\partial
B_\sigma(0)}|y_N|^{2\alpha}u\left(\frac{y}{|y|}\right)\left((\nabla_x
u)\left(\frac{y}{|y|}\right)\cdot \frac{y}{|y|}\right)\nonumber\\
&=&\frac{1}{\sigma^{2N-2+4\alpha}}\int_{\partial
B_{1/\sigma}(0)}\frac{|x_N|^{2\alpha}}{|x|^{4\alpha}}u(x)\left((\nabla_x
u)\left(x\right)\cdot \frac{x}{|x|}\right)\sigma^{2(N-1)}\nonumber\\
&=&\int_{\partial
B_{1/\sigma}(0)}|x_N|^{2\alpha}u(x)\left((\nabla_x
u)\left(x\right)\cdot \frac{x}{|x|}\right)\nonumber\\
&=&\int_{\partial B_{1/\sigma}(0)}|x_{N}|^{2\alpha}\cdot
u\frac{\partial u}{\partial\mathbf{n}}.
\end{eqnarray}
As above, we deduce that
$$\lim_{\sigma\rightarrow 0}\int_{\partial B_\sigma(0)}|y_N|^{2\alpha}\widetilde{u}\nabla_y
\widetilde{u}\cdot\frac{y}{|y|}=0,\ \lim_{\sigma\rightarrow
0}(N-2+2\alpha)\int_{\partial
B_\sigma(0)}\frac{|y_N|^{2\alpha}}{|y|} \widetilde{u}^2(y)=0.$$
Then by (\ref{ggfbjjfufhhhf77qe}) and (\ref{iidhdhgtdrddfd}), we
get that
\begin{eqnarray}\label{hhfbvggyfg99id}
\lim_{\sigma\rightarrow 0}\int_{\partial
B_{1/\sigma}(0)}|x_{N}|^{2\alpha}\cdot u\frac{\partial
u}{\partial\mathbf{n}}=0.
\end{eqnarray}
By (\ref{hfhbcgdteyye00}), (\ref{hfb77fyfgggg}) and
(\ref{hhfbvggyfg99id}), we get that
\begin{eqnarray}\label{jjjjdvggcgd}
\lim_{\rho\rightarrow \infty}\int_{\partial
B_{\rho}}B(\rho,x,u,\nabla u)=0.
\end{eqnarray}
By the fact that $\widetilde{u}\in L^\infty(\mathbb{R}^N)$, we
deduce that $u$ satisfies $|x|^{-(N-2+2\alpha)}$ decay at
infinity. Thus,
$$\lim_{\rho\rightarrow\infty}\int_{\partial
B_{\rho}(0)}(x\cdot\mathbf{n})\cdot K(x)
|x_{N}|^{\alpha\cdot2^{*}(s)-s}|u|^{2^{*}(s)}=0.$$ Then by
(\ref{jjjjdvggcgd}) and (\ref{xzc4423}), we obtain the desired
result of this Lemma.\hfill$\Box$

\begin{remark}\label{ygfhcvvgf8y}
Suppose  $K, \widehat{K}\in C^{1}(\mathbb{R}^{N})\cap
L^{\infty}(\mathbb{R}^{N})$. If $x\nabla K(x)\geq 0$ or $x\nabla
K(x)\leq 0$, $\forall x\in \mathbb{R}^{N}$ and $x\nabla
K(x)\not\equiv0,$ then equation (\ref{yyyhsf}) does not have
solution lying  in $ X_{\alpha }(\mathbb{R}^{N}).$ Furthermore, by
Remark \ref{ccv5} and this corollary, we know that if $x\nabla
K(x)\geq 0$ or $x\nabla K(x)\leq 0$, $\forall x\in \mathbb{R}^{N}$
and $x\nabla K(x)\not\equiv0,$ then equation (\ref{yyyhsf}) does
not have positive solution  in $ X_{\alpha, loc
}(\mathbb{R}^{N}).$
\end{remark}

\begin{Theorem}\label{cet551}
Suppose that $K\in C^{1}(\mathbb{R}^{N})$. If $u\in
X_{\alpha}(\mathbb{R}^{N})$ is a weak solution of equation
(\ref{yyyhsf}),  then
$$\int_{\mathbb{R}^{N}}\frac{\partial K}{\partial x_{i}}\cdot
|x_{N}|^{\alpha\cdot 2^{*}(s)-s} |u|^{2^{*}(s)}=0,\ 1\leq i\leq
N-1.$$
\end{Theorem}
\noindent{\bf Proof.} Let $\varphi_{R}$ be a cut-off function
which satisfies that $\varphi_{R}\equiv 1$ in $B_{R}(0)$,
$\varphi_{R}\equiv 0$ in $\mathbb{R}^{N}\setminus B_{R+1}(0)$ and
$|\nabla\varphi_{R}(x)|\leq 1,$ $\forall x\in\mathbb{R}^{N}.$ For
$1\leq i\leq N-1,$ multiplying the equation (\ref{yyyhsf}) by
$\varphi_{R}\frac{\partial u}{\partial x_{i}}$ and integrating in
$\{x_{N}>\epsilon\}$, we obtain
\begin{eqnarray}\label{bcvtre5}
-\int_{x_{N}>\epsilon}div(|x_{N}|^{2\alpha}\nabla
u)\cdot\left(\varphi_{R}\frac{\partial u}{\partial x_{i}}\right)
=\int_{x_{N}>\epsilon}K(x)|x_{N}|^{\alpha\cdot
2^{*}(s)-s}|u|^{2^{*}(s)-2}u\cdot\left(\varphi_{R}\frac{\partial
u}{\partial x_{i}}\right).
\end{eqnarray}
Through integrating by parts, we get
\begin{eqnarray}\label{vcfer1}
&&-\int_{x_{N}>\epsilon}div(|x_{N}|^{2\alpha}\nabla
u)\cdot\left(\varphi_{R}\frac{\partial u}{\partial
x_{i}}\right)\nonumber\\
&=&\int_{x_{N}=\epsilon}|x_{N}|^{2\alpha}\varphi_{R}\frac{\partial
u}{\partial x_{i}}\cdot\frac{\partial u}{\partial x_{N}}
+\int_{x_{N}>\epsilon}|x_{N}|^{2\alpha}\nabla u\cdot\nabla
\left(\varphi_{R}\frac{\partial u}{\partial
x_{i}}\right)\nonumber\\
&=&\int_{x_{N}=\epsilon}|x_{N}|^{2\alpha}\varphi_{R}\frac{\partial
u}{\partial x_{i}}\cdot\frac{\partial u}{\partial x_{N}}+
\int_{x_{N}>\epsilon}|x_{N}|^{2\alpha}\varphi_{R}\cdot\nabla
u\cdot\nabla \left(\frac{\partial u}{\partial
x_{i}}\right)\nonumber\\
&&+\int_{x_{N}>\epsilon}|x_{N}|^{2\alpha}\frac{\partial
u}{\partial x_{i}}\cdot\nabla u\cdot\nabla\varphi_{R}.
\end{eqnarray}
Since
\begin{eqnarray}\label{nnn7y}
&&\int_{x_{N}>\epsilon}|x_{N}|^{2\alpha}\varphi_{R}\cdot\nabla
u\cdot\nabla \left(\frac{\partial u}{\partial
x_{i}}\right)\nonumber\\
&=&\sum^{N}_{j=1}\int_{x_{N}>\epsilon}|x_{N}|^{2\alpha}\varphi_{R}\cdot\frac{\partial
u}{\partial x_{j}}\frac{\partial^{2}u}{\partial x_{i}\partial
x_{j}}\nonumber\\
&=&-\sum^{N}_{j=1}\int_{x_{N}>\epsilon}\frac{\partial u}{\partial
x_{j}}\cdot\frac{\partial}{\partial
x_{i}}\left(|x_{N}|^{2\alpha}\varphi_{R}\frac{\partial u}{\partial
x_{j}}\right) \quad (\mbox{though integral by parts})\nonumber\\
&=&-\sum^{N}_{j=1}\int_{x_{N}>\epsilon}\frac{\partial u}{\partial
x_{j}}\cdot|x_{N}|^{2\alpha}\varphi_{R}\frac{\partial^{2}u}{\partial
x_{i}\partial
x_{j}}-\sum^{N}_{j=1}\int_{x_{N}>\epsilon}\left(\frac{\partial
u}{\partial x_{j}}\right)^{2}\cdot
|x_{N}|^{2\alpha}\cdot\frac{\partial\varphi_{R}}{\partial
x_{i}}\nonumber\\
&=&-\int_{x_{N}>\epsilon}|x_{N}|^{2\alpha}\varphi_{R}\cdot\nabla
u\cdot\nabla \left(\frac{\partial u}{\partial
x_{i}}\right)-\int_{x_{N}>\epsilon}|x_{N}|^{2\alpha}|\nabla
u|^{2}\cdot\frac{\partial \varphi_{R}}{\partial x_{i}},\nonumber
\end{eqnarray}
we get that
\begin{equation}\label{cxcsde}
\int_{x_{N}>\epsilon}|x_{N}|^{2\alpha}\varphi_{R}\cdot\nabla
u\cdot\nabla \left(\frac{\partial u}{\partial x_{i}}\right)
=-\frac{1}{2}\int_{x_{N}>\epsilon}|x_{N}|^{2\alpha}|\nabla
u|^{2}\cdot\frac{\partial \varphi_{R}}{\partial x_{i}}.
\end{equation}
By (\ref{vcfer1}) and (\ref{cxcsde}), we get that
\begin{eqnarray}\label{bcvfd5r5r}
&&-\int_{x_{N}>\epsilon}div(|x_{N}|^{2\alpha}\nabla
u)\cdot\left(\varphi_{R}\frac{\partial u}{\partial
x_{i}}\right)\nonumber\\
&=&\int_{x_{N}=\epsilon}|x_{N}|^{2\alpha}\varphi_{R}\frac{\partial
u}{\partial x_{i}}\cdot\frac{\partial u}{\partial
x_{N}}-\frac{1}{2}\int_{x_{N}>\epsilon}|x_{N}|^{2\alpha}|\nabla
u|^{2}\cdot\frac{\partial \varphi_{R}}{\partial
x_{i}}+\int_{x_{N}>\epsilon}|x_{N}|^{2\alpha}\frac{\partial
u}{\partial x_{i}}\cdot\nabla u\cdot\nabla\varphi_{R}.\nonumber\\
\end{eqnarray}
In a similar manner, we have
\begin{eqnarray}\label{bcv4re}
&&-\int_{x_{N}<-\epsilon}div(|x_{N}|^{2\alpha}\nabla
u)\cdot\left(\varphi_{R}\frac{\partial u}{\partial
x_{i}}\right)\nonumber\\
&=&-\int_{x_{N}=-\epsilon}|x_{N}|^{2\alpha}\varphi_{R}\frac{\partial
u}{\partial x_{i}}\cdot\frac{\partial u}{\partial
x_{N}}-\frac{1}{2}\int_{x_{N}<-\epsilon}|x_{N}|^{2\alpha}|\nabla
u|^{2}\cdot\frac{\partial \varphi_{R}}{\partial
x_{i}}\nonumber\\
&&+\int_{x_{N}<-\epsilon}|x_{N}|^{2\alpha}\frac{\partial
u}{\partial x_{i}}\cdot\nabla u\cdot\nabla\varphi_{R}.
\end{eqnarray}
Adding (\ref{bcvfd5r5r}) and (\ref{bcv4re}) and using the fact
that $ \lim_{\epsilon\rightarrow
0}\int_{x_{N}=\pm\epsilon}|x_{N}|^{2\alpha}\varphi_{R}\frac{\partial
u}{\partial x_{i}}\cdot\frac{\partial u}{\partial x_{N}}=0$ (
 by  Proposition \ref{tqaaa3w} and  \ref{fff1fc}), we get
that ( letting $\epsilon\rightarrow 0$ ),
\begin{eqnarray}\label{vcvrtw43}
&&-\int_{\mathbb{R}^{N}}div(|x_{N}|^{2\alpha}\nabla
u)\cdot\left(\varphi_{R}\frac{\partial u}{\partial
x_{i}}\right)\nonumber\\
&=&-\frac{1}{2}\int_{\mathbb{R}^{N}}|x_{N}|^{2\alpha}|\nabla
u|^{2}\cdot\frac{\partial \varphi_{R}}{\partial
x_{i}}+\int_{\mathbb{R}^{N}}|x_{N}|^{2\alpha}\frac{\partial
u}{\partial x_{i}}\cdot\nabla u\cdot\nabla\varphi_{R}.
\end{eqnarray}
By (\ref{vcvrtw43}), we have
\begin{equation}\label{v5r1}
\lim_{R\rightarrow\infty}\int_{\mathbb{R}^{N}}div(|x_{N}|^{2\alpha}\nabla
u)\cdot\left(\varphi_{R}\frac{\partial u}{\partial
x_{i}}\right)=0.
\end{equation}
On the other hand,
\begin{eqnarray}\label{xcde1}
&&\int_{\mathbb{R}^{N}}K(x)|x_{N}|^{\alpha\cdot
2^{*}(s)-s}|u|^{2^{*}(s)-2}u\cdot\left(\varphi_{R}\frac{\partial
u}{\partial x_{i}}\right)\nonumber\\
&=&\frac{1}{2^{*}(s)}\int_{\mathbb{R}^{N}}K(x)|x_{N}|^{\alpha\cdot
2^{*}(s)-s}\varphi_{R}\frac{\partial}{\partial
x_{i}}\left(|u|^{2^{*}(s)}\right)\nonumber\\
 &=&-\frac{1}{2^{*}(s)}\int_{\mathbb{R}^{N}}|u|^{2^{*}(s)}\cdot\frac{\partial}{\partial
x_{i}}\left(K(x)|x_{N}|^{\alpha\cdot 2^{*}(s)-s}
\varphi_{R}\right)\nonumber\\
&=&-\frac{1}{2^{*}(s)}\int_{\mathbb{R}^{N}}\frac{\partial
K}{\partial x_{i}} |x_{N}|^{\alpha\cdot
2^{*}(s)-s}|u|^{2^{*}(s)}\varphi_{R}-
\frac{1}{2^{*}(s)}\int_{\mathbb{R}^{N}}K(x) |x_{N}|^{\alpha\cdot
2^{*}(s)-s}|u|^{2^{*}(s)}\cdot\frac{\partial\varphi_{R}}{\partial
x_{i}}.\nonumber
\end{eqnarray}
Letting $R\rightarrow\infty$ in the above identity, we obtain
\begin{equation}\label{9910}
\lim_{R\rightarrow\infty}\int_{\mathbb{R}^{N}}K(x)|x_{N}|^{\alpha\cdot
2^{*}(s)-s}|u|^{2^{*}(s)-1}\cdot\left(\varphi_{R}\frac{\partial
u}{\partial
x_{i}}\right)=-\frac{1}{2^{*}(s)}\int_{\mathbb{R}^{N}}\frac{\partial
K}{\partial x_{i}} |x_{N}|^{\alpha\cdot 2^{*}(s)-s}|u|^{2^{*}(s)}.
\end{equation}
By (\ref{v5r1}), (\ref{9910}) and (\ref{bcvtre5}), we get
$\int_{\mathbb{R}^{N}}\frac{\partial K}{\partial x_{i}}
|x_{N}|^{\alpha\cdot 2^{*}(s)-s}|u|^{2^{*}(s)}=0.$ \hfill$\Box$

\medskip

The following kind of result will be used in the blow-up analysis
of equation (\ref{yyyhsf}). Similar results have been used in
\cite{Liyy}.

\begin{proposition}\label{byttt}
\begin{itemize}
\item[(i).] For $u=|x|^{-(N-2+2\alpha)},$ $B(\sigma,x,u,\nabla
u)=0$ for all $x\in\partial B_{\sigma}(0)$; \item[(ii).] For
$u(x)=|x|^{-(N-2+2\alpha)}+A+\xi(x),$ with $A>0$ and $\xi(0)=0,$
$-div(|x_{N}|^{2\alpha}\nabla\xi)=0$ weakly in $B_{1}(0),$
 there exists $\overline{\sigma}$ such that
$$B(\sigma,x,u,\nabla
u)<0\ \mbox{for all} \ x\in\partial B_{\sigma}(0)\ \mbox{and}\
0<\sigma<\overline{\sigma}$$ and
$$\lim_{\sigma\rightarrow 0}\int_{\partial B_{\sigma}(0)}B(\sigma,x,u,\nabla
u)=-\frac{1}{2}A(N-2+2\alpha)^{2}\int_{\partial
B_{1}(0)}|x_{N}|^{2\alpha}.$$ \end{itemize}
\end{proposition}
\noindent{\bf Proof.} (i). By straightforward calculation, we have
$\nabla u=-(N-2+2\alpha)|x|^{-(N-2+2\alpha)-2}x$ and
$\frac{\partial u}{\partial \mathbf{n}}=\mathbf{n}\cdot\nabla
u=\frac{x}{|x|}\cdot\nabla
u=-(N-2+2\alpha)|x|^{-(N-2+2\alpha)-1}.$ It follows that $|\nabla
u|^{2}|_{|x|=\sigma}=(N-2+2\alpha)^{2}\sigma^{-2(N-2+2\alpha)-2},$
$u\frac{\partial u}{\partial
\mathbf{n}}|_{|x|=\sigma}=-(N-2+2\alpha)\sigma^{-2(N-2+2\alpha)-1}$
and $(\partial
u/\partial\mathbf{n})^{2}|_{|x|=\sigma}=(N-2+2\alpha)^{2}\sigma^{-2(N-2+2\alpha)-2}.$
Thus $B(\sigma,x,u,\nabla u)=0$ for all $x\in\partial
B_{\sigma}(0)$.

\medskip

(ii). From the assumptions $\xi$ holds and the regularity result
in section \ref{tgrfe} (see Proposition \ref{tqaaa3w}), we know
that $\xi$ and $\frac{\partial\xi}{\partial x_{i}}$ ($1\leq i\leq
N-1$) are local H\"older continuous in $B_{1}(0)$. And $\xi$ is
$C^{2}$ continuous in $B_{1}(0)\setminus\{x_{N}=0\}$.
Straightforward calculation shows that $\nabla
u(x)=-(N-2+2\alpha)|x|^{-(N-2+2\alpha)-2}x+\nabla\xi(x)$ and
\begin{eqnarray}\label{cz88u}
 &&\frac{\partial u}{\partial\mathbf{n}}=\frac{x}{|x|}\cdot\nabla
u
=-(N-2+2\alpha)|x|^{-(N-2+2\alpha)-1}+\frac{x}{|x|}\cdot\nabla\xi(x).\nonumber
\end{eqnarray}
Then by straightforward calculation and using the result  $(i)$ of
this proposition, we get
\begin{eqnarray}\label{bcv556a}
 &&B(\sigma,x,u,\nabla u)\Big|_{|x|=\sigma}\nonumber\\
&=& B(\sigma,x,|x|^{-(N-2+2\alpha)},\nabla
(|x|^{-(N-2+2\alpha)}))\Big|_{|x|=\sigma}+B(\sigma,x,\xi,\nabla
\xi)\Big|_{|x|=\sigma}\nonumber\\
&&-\frac{(N-2+2\alpha)^{2}}{2}A\sigma^{-(N-2+2\alpha)-1}|x_{N}|^{2\alpha}+R_{\sigma}\nonumber\\
&=&B(\sigma,x,\xi,\nabla \xi)\Big|_{|x|=\sigma}
-\frac{(N-2+2\alpha)^{2}}{2}A\sigma^{-(N-2+2\alpha)-1}|x_{N}|^{2\alpha}+R_{\sigma},
\end{eqnarray}
where $R_{\sigma}$ equals to \begin{eqnarray}\label{ttegdfdhhh}
&&|x_{N}|^{2\alpha}(
\frac{N-2+2\alpha}{2}A\sigma^{-1}(x\cdot\nabla\xi(x))
-\frac{(N-2+2\alpha)^{2}}{2}\sigma^{-(N-2+2\alpha)-1}\xi(x)\nonumber\\
&&\quad\quad\quad\quad
-\frac{N-2+2\alpha}{2}\sigma^{-(N-2+2\alpha)-1} \cdot
(x\cdot\nabla\xi(x)))\Big|_{|x|=\sigma}.\nonumber
\end{eqnarray}
By the regularity results   of $\xi$ and the condition $\xi(0)=0$,
we deduce that
\begin{eqnarray}\label{yyetdff7y00oi}B(\sigma,x,u,\nabla u)\Big|_{|x|=\sigma}<0\
\mbox{if}\ \sigma\ \mbox{small enough.}
\end{eqnarray}
Multiplying equation $-div(|x_{N}|^{2\alpha}\nabla \xi)=0$ by $1$
and integrating in $B_{\sigma}(0),$ we have
\begin{eqnarray}\label{vcrr43}
0=-\int_{B_{\sigma}(0)}div(|x_{N}|^{2\alpha}\nabla\xi)\cdot
1=-\int_{\partial
B_{\sigma}(0)}|x_{N}|^{2\alpha}\frac{\partial\xi}{\partial\mathbf{n}}=-\sigma^{-1}\int_{\partial
B_{\sigma}(0)}|x_{N}|^{2\alpha}(x\cdot\nabla\xi).
\end{eqnarray}
Moreover, by $\xi(0)=0,$ we can get that
\begin{eqnarray}\label{cxvvfdre}
\lim_{\sigma\rightarrow 0}\sigma^{-(N-2+2\alpha)-1}\int_{\partial
B_{\sigma}(0)}|x_{N}|^{2\alpha}\xi(x)=0.
\end{eqnarray}
By (\ref{vcrr43}) and (\ref{cxvvfdre}), we get that
\begin{eqnarray}\label{99iqjsvvx} \lim_{\sigma\rightarrow
0}\int_{\partial B_{\sigma}(0)}R_{\sigma}=0.
\end{eqnarray}
Let $K\equiv 0$ in Theorem \ref{bvcxq}, we obtain
\begin{equation}\label{iudt}
\int_{\partial B_{\sigma}(0)}B(\sigma,x,\xi,\nabla \xi)=0.
\end{equation}
 Thus by
  $(\ref{99iqjsvvx})$ and (\ref{iudt}),
we get that
\begin{eqnarray}\label{v2}
\lim_{\sigma\rightarrow 0}\int_{\partial
B_{\sigma}(0)}B(\sigma,x,u,\nabla
u)&=&-\frac{1}{2}A(N-2+2\alpha)^{2}\lim_{\sigma\rightarrow
0}\sigma^{-(N-2+2\alpha)-1}\int_{\partial
B_{\sigma}(0)}|x_{N}|^{2\alpha}\nonumber\\
&=&-\frac{1}{2}A(N-2+2\alpha)^{2}\int_{\partial
B_{1}(0)}|x_{N}|^{2\alpha}<0.
\end{eqnarray}
The result of this Proposition follows from (\ref{yyetdff7y00oi})
and (\ref{v2}).\hfill$\Box$

\bigskip

{\noindent\bf Acknowledgement.} This work was supported by the
Natural Science Foundation of China (No. 10901112).

\end{document}